\documentclass[oneside,leqno,11pt]{article}
\usepackage{amssymb,amsmath,latexsym}

\def\ga{\mathfrak{a}}

\def\ge{\mathfrak{e}}

\def\gg{\mathfrak{g}}
\def\gh{\mathfrak{h}}

\def\gk{\mathfrak{k}}
\def\gl{\mathfrak{l}}

\def\gp{\mathfrak{p}}
\def\gq{\mathfrak{q}}
\def\gr{\mathfrak{r}}
\def\gs{\mathfrak{s}}

\def\ggg{> \hskip -5 pt >}

\def\C{\mathbb{C}}

\def\E{\mathbb{E}}

\def\H{\mathbb{H}}

\def\K{\mathbb{K}}
\def\L{\mathbb{L}}

\def\P{\mathbb{P}}

\def\R{\mathbb{R}}

\def\Z{\mathbb{Z}}

\def\cA{\mathcal{A}}
\def\cB{\mathcal{B}}
\def\cC{\mathcal{C}}

\def\cE{\mathcal{E}}
\def\cF{\mathcal{F}}

\def\cI{\mathcal{I}}

\def\cL{\mathcal{L}}

\def\cO{\mathcal{O}}

\def\cR{\mathcal{R}}
\def\cS{\mathcal{S}}

\newtheorem{theorem}[equation]{Theorem}

\newtheorem{lemma}[equation]{Lemma}

\newtheorem{corollary}[equation]{Corollary}

\newtheorem{proposition}[equation]{Proposition}

\newtheorem{definition}[equation]{Definition}

\newtheorem{example}[equation]{Example}

\newtheorem{remark}[equation]{Remark}

\newtheorem{summary}[equation]{Summary}

\title{Cycle Spaces of Flag Domains: A Complex Geometric Viewpoint}

\author{Alan T. Huckleberry
\thanks{
Research partially supported by
Schwerpunkt ``Global methods in complex
geometry'' and SFB-237 of the Deutsche Forschungsgemeinschaft.} \; \& \;
Joseph A. Wolf
\thanks{
Research partially supported by NSF Grant
DMS 99-88643 and by the SFB-237 of the Deutsche Forschungsgemeinschaft.}}

\date{29 October 2002}

\begin{document}

\maketitle

\abstract{This is a survey of history, methods and 
developments in the theory of cycle spaces of flag domains, and
new results on double fibration transforms and their applications.
}

{\scriptsize
$$
\begin{aligned}
\text{\large \bf Contents:}& \\
\S 0. &\text{ Introduction} \\
\text{Part I:} &\text{ Background} \\
\S 1. &\text{ Flag Domains and Compact Subvarieties} \\
   &\phantom{abcd} \text{Measurability} \\
   &\phantom{abcd} \text{Compact subvarieties} \\
\S 2. &\text{ Basic Facts on the Cycle Space} \\
   &\phantom{abcd} \text{Hermitian trichotomy} \\
\S 3. &\text{ The Exhaustion Function and } (q+1)\text{--Completeness} \\
   &\phantom{abcd} \text{Cohomology vanishing theorems} \\
   &\phantom{abcd} \text{The Stein property for cycle spaces of measurable open orbits} \\
\S 4. &\text{ Early Problems and Results on the Double Fibration Transform} \\
   &\phantom{abcd} \text{Double fibration} \\
   &\phantom{abcd} \text{Pull--back} \\
   &\phantom{abcd} \text{Push--down} \\
   &\phantom{abcd} \text{Flag domain case} \\
\text{Part II:} &\text{ The Complex Geometric Approach} \\
\S 5. &\text{ Introduction to the Complex Geometric Approach} \\
\S 6. &\text{ The Equivalences $\Omega _{adpt}=\Omega _{AG}=\Omega _I$} \\
   &\phantom{abcd} \text{Adapted complex structures} \\
   &\phantom{abcd} \text{Basic properties of plurisubharmonic functions} \\
   &\phantom{abcd} \text{The adapted structure for Riemannian symmetric spaces} \\
   &\phantom{abcd} \text{Proper actions} \\
   &\phantom{abcd} \text{Incidence geometry and the domain $\Omega _I$} \\
   &\phantom{abcd} \text{Domains defined by invariant hypersurfaces} \\
   &\phantom{abcd} \text{The equality $\Omega _I=\Omega _{AG}$} \\
\end{aligned}
$$

$$
\begin{aligned}
\text{\large \bf Contents:}&\text{ (continued)} \\
\S 7. &\text{ Transversal Schubert Varieties} \\
   &\phantom{abcd} \text{Duality} \\
   &\phantom{abcd} \text{Triality} \\
   &\phantom{abcd} \text{The equality $\Omega _W(D)=\Omega _S(D)$} \\
   &\phantom{abcd} \text{Spaces of cycles in lower dimensional $G_0$--orbits} \\
\S 8. &\text{ Cycle Domains in the Hermitian Case} \\
\S 9. &\text{ Kobayashi Hyperbolicity} \\
   &\phantom{abcd} \text{Families of hypersurfaces} \\
   &\phantom{abcd} \text{Invariant hyperbolic domains in $\Omega $} \\
\S 10. &\text{ The Maximal Domain of Hyperbolicity} \\
   &\phantom{abcd} \text{The linear model} \\
   &\phantom{abcd} \text{Orbit structure} \\
   &\phantom{abcd} \text{Existence of a $Q_2$-slice} \\
   &\phantom{abcd} \text{Analysis of a $Q_2$-slice} \\
   &\phantom{abcd} \text{Characterization of cycle domains} \\
\text{Part III:} &\text{ Applications and Open Problems} \\
\S 11. &\text{ Recent Results on the Double Fibration Transform} \\
\S 12. &\text{ Unitary Representations of Real Reductive Lie Groups} \\
\S 13. &\text{ Variation of Hodge Structure}
\end{aligned}
$$
}

\setcounter{section}{-1}
\section{Introduction} \label{sec0}
\setcounter{equation}{0}

Cycle space theory is a basic chapter in complex analysis.  Since
since the 1960's its importance has been underlined by its  role in the geometry of
flag domains and applications, by means of double fibration transforms, to
variation of Hodge structure and to the representation theory of semisimple
Lie groups.  This developed very slowly until a few of years ago
when methods of complex analytic geometry, in particular the methods of
Schubert slices, Schubert domains, Iwasawa domains and supporting 
hypersurfaces, were introduced.  Early in 2002
those methods were used to settle a number of outstanding questions.
This effectively enabled the use of double fibration transforms in all flag
domain situations.  This has very interesting consequences for geometric
construction of representations of semisimple Lie groups, especially for the
construction of singular representations.  It also has many potential
interesting consequences for automorphic cohomology and other aspects of
variation of Hodge structure.  In this article we survey the recent results,
filling in the background as necessary, and present some new results that 
help to complete the picture. 

Part I,  ``Background'', is an exposition of
flag domains and their cycle spaces before the introduction of the new
complex geometric methods.  Section \ref{sec1} recalls the very basic
results on flag domains and compact subvarieties.  Section \ref{sec2} goes
into the complex structure of these cycle spaces and describes the three
basic possibilities.  Section \ref{sec3} describes a particular exhaustion
function for measurable flag domains, and its consequences for cohomology
vanishing theorems that are crucial to the double fibration transform and
application to semisimple representation theory.  The basic aspects of that 
double fibration transform are described in Section \ref{sec4}.

Part II, ``The Complex Geometric Approach'', introduces the methods of
Schubert slices, Schubert domains, Iwasawa domains and supporting
hypersurfaces, and the use of Kobayashi hyperbolicity in this context.  
In order to orient the reader who is not working in complex
analysis, Section \ref{sec5} is an introduction to the methods and ideas, and a
sketch of the remainder of Part II.  Section \ref{sec6} introduces
several related domains, and proves certain equivalences among them,
specifically $\Omega_{adpt} \cong \Omega_{AG} \cong  \Omega_I$.  Section
\ref{sec7} goes into the key notion of transversal Schubert varieties, 
explaining an enhanced duality theory and the argument that the cycle
space $\Omega_W(D)$ is equal to the Schubert domain $\Omega_S(D)$ in all but
certain exceptional cases related to hermitian symmetric spaces.  The situation
of hermitian symmetric spaces is completely described in Section \ref{sec8}.
Another new element in this picture, the use of Kobayashi hyperbolicity, is
described in Section \ref{sec9}.  In Section \ref{sec10} these ingredients
are combined to describe the maximal domain of hyperbolicity; the result
is $\Omega_{AG} \cong  \Omega_I \cong \Omega_D(D) \cong \Omega_W(D)$,
again except in certain exceptional cases related to hermitian symmetric spaces
which are described completely in  Section \ref{sec8}.

Part III, ``Applications and Open Problems'', applies the results of Part II
to the mechanism of the double fibration transform, and discusses certain
applications.  The material on the double fibration transform, key parts of
which are new, appears in Section \ref{sec11}.  Consequences for representations
of real reductive Lie groups are discussed in Section \ref{sec12}, and in
Section \ref{sec13} there is a discussion of variation of Hodge structure and
automorphic cohomology.

\bigskip
\centerline{\Large \bf \underline{Part I: Background.}}
\bigskip

In this Part we describe the early results on the cycle space and the
double fibration transform.  For the most part those results are based on
Lie structure theory.

\section{Flag Domains and Compact Subvarieties.} \label{sec1}
\setcounter{equation}{0}

We begin by reviewing the 
basic setup for flag domains and their maximal compact subvarieties
which was presented in Wolf \cite{W2}.  We review the part of 
\cite{W2} that is relevant to the theory of cycle spaces of flag domains.
\medskip

Let $G$ be a complex semisimple Lie group and $Q$ a parabolic subgroup. 
The compact algebraic homogeneous space 
$Z = G/Q$ is called a complex flag manifold.    
Write $\gg$ and $\gq$ for the respective Lie algebras of
$G$ and $Q$.  Then $Q$ is the $G$--normalizer of $\gq$.  Thus we may
view $Z$ as the set of $G$--conjugates of $\gq$.  The correspondence
is $z \leftrightarrow \gq_z$ where $\gq_z$ is the Lie algebra
of the isotropy subgroup $Q_z$ of $G$ at $z$. 

Let $G_0$ be a real form of $G$ in the sense that 
there is a homomorphism $\varphi : G_0 \to G$ 
such that $\varphi(G_0)$ is closed in $G$ and $d\varphi : \gg_0 \to \gg$ is an
isomorphism onto a real form of $\gg$.  In this paper we will only consider the situation
where $\varphi$ is an inclusion, $\varphi : G_0 \hookrightarrow G$, so
we now assume $G_0 \subset G$ and that $G_0$ is noncompact.

Write $g \mapsto \overline{g}$ for
complex conjugation of $G$ over $G_0$ and of $\gg$ over 
$\gg_0$\,.  We recall some of the basic
facts about $G_0$--orbits on Z. 

If $z \in Z$ then $\gq_z \cap \overline{\gq_z}$ contains a Cartan
subalgebra $\gh$ of $\gg$.  We may assume that $\gh =
\overline{\gh}$, in other words that $\gh$ is the complexification
of a Cartan subalgebra $\gh_0 = \gh \cap \gg_0$ of 
$\gg_0$\,.  There is a choice of positive root system
$\Delta^+ = \Delta^+(\gg, \gh)$ such that $\gq_z$ is
the standard parabolic subalgebra $\gq_\Phi$ defined by some subset 
$\Phi \subset \Psi$ where $\Psi = \Psi(\gg, \gh, \Delta^+)$ is 
the corresponding simple root system.  In other words, $\gq_z = 
\gq_\Phi$ where
\begin{equation}
\begin{aligned}
{}& \Phi^r = \{\alpha \in \Delta \mid \alpha \text{ is a linear combination of 
 	elements of } \Phi \}, \\
{}& \Phi^n = \{\alpha \in \Sigma^+ \mid \alpha \notin \Phi^r\}, \text{ and } \\
& \gq_\Phi = \gq_\Phi^r + \gq_\Phi^{-n} \text{ with }
 \gq_\Phi^r = \gh + \sum_{\alpha \in \Phi^r} \gg_\alpha \text{ and }
 \gq_\Phi^{-n} = \sum_{\alpha \in \Phi^n} \gg_{-\alpha} \ . 
\end{aligned}
\end{equation}
It follows that $G_0$ acts on $Z$ with
only finitely many orbits; in particular there are open orbits.
We refer to the open orbits as
{\sl flag domains}.  As $G_0$--invariant open subsets of $Z$, the flag 
domains $D \subset Z$ are $G_0$--homogeneous complex manifolds.

\noindent {\bf Measurability.}
\medskip

A flag domain $D = G_0(z) \subset Z$ is called {\sl measurable} if it 
carries a $G_0$--invariant volume element.  This is the type of flag
domain currently of most interest in representation theory.  
More precisely, the following conditions are equivalent:

\noindent (1.2a) \hskip .2 cm The orbit $G_0(z)$ is measurable. \hfill\newline
\noindent (1.2b) \hskip .2 cm $G_0\cap Q_z$ is the $G_0$--centralizer of a (compact)
torus subgroup of $G_0$\,. \hfill\newline
\noindent (1.2c) \hskip .2 cm $D$ has a $G_0$--invariant, possibly--indefinite, 
K\" ahler metric, thus a $G_0$--invariant \hfill\newline
\indent \indent \indent measure obtained from the volume form of 
that metric. \hfill\newline
\noindent (1.2d) \hskip .2 cm $\overline{\Phi^r} = \Phi^r$, and 
$\overline{\Phi^n} = - \Phi^n$
where $\gq_z = \gq_\Phi$\,. \hfill\newline
\noindent (1.2e) \hskip .2 cm $\gq_z \cap \overline{\gq_z}$ is 
reductive, i.e. $\gq_z \cap \overline{\gq_z} = 
\gq_z^r \cap \overline{q_z^r}$\,.  \hfill\newline
\noindent (1.2f) \hskip .2 cm $\gq_z \cap \overline{\gq_z}
= \gq_z^r$ \,.\hfill\newline
\noindent (1.2g) \hskip .2 cm $\overline{\gq}$ is {\rm Ad\,}(G)--conjugate to the
parabolic subalgebra $\gq^r + \gq^n$ opposite to
$\gq$. \hfill\newline
\setcounter{equation}{2}
\noindent In particular, since (1.2g) is independent of choice of $z$,
if one open $G_0$--orbit on $Z$ is measurable
then all open $G_0$--orbits are measurable.

Condition (1.2d) holds whenever the Cartan subalgebra
$\gh_0 = \gh \cap \gg_0$ of $\gg_0$ 
corresponds to a compact Cartan subgroup $H_0 \subset G_0$\,.  
(Here $\gh = \overline{\gh}$ is the Cartan subalgebra relative to 
which $\gq_z = \gq_\Phi$\,.)  For in that case $\overline{\alpha}
= -\alpha$ for every $\alpha \in \Delta(\gg, \gh)$.  In particular,
if $G_0$ has discrete series representations (so that by a result of
Harish--Chandra it has a compact Cartan subgroup) then every open
$G_0$--orbit on $Z$ is measurable.
Condition (1.2d) is also automatic if $Q$ is a Borel subgroup of $G$, and 
more generally Condition (1.2g) provides a quick test for measurability.

\noindent {\bf Compact subvarieties.} 
\medskip

We now fix $z \in Z$ such that $D = G_0(z)$ is open in $Z$.
For convenience we suppose that $z$ is the base point in $Z = G/Q$, so 
$Q = Q_z$ and $\gq = \gq_z$\,.  For notational consistency with 
many papers in this area, we write $L$ for the Levy
component $Q^r$ of $Q$.  
So $D$ is measurable if and only if $Q \cap G_0$ is a real form $L_0$
of $L$, and in that case $D \cong G_0/L_0$\,.  

Fix a Cartan involution $\theta$ of $G_0$ that stabilizes the Cartan subgroup
$H_0 \subset G_0$\,, and denote its fixed point sets on $G_0$ and $G$ by
$K_0 = G_0^\theta$ and $K = G^\theta$.  Then $K_0$ is a maximal compact
subgroup of $G_0$ and $K$ is its complexification.   $L \cap K_0$ is a real 
form of $L \cap K$ and $K_0(z) \cong K_0/(L \cap K_0)$.

As $D$ is open we may 
assume $\gh$ chosen so that $H_0 \cap K_0$ is a Cartan subgroup of
$K_0$\,, in other words so that $H_0$ is a {\sl fundamental} Cartan
subgroup of $G_0$\,.   Use $\gh$ for the standard Weyl basis construction
of a $\theta$--stable compact real form $\gg_u \subset \gg$.  Then
$G_0 \cap G_u = K_0$ and $\gk = (\gk \cap \gl) + 
(\gk \cap \gr_-) + (\gk \cap \gr_+)$.  Thus 
$K(z) \cong K/(K \cap Q)$ is a complex flag submanifold of $Z$, and $K_0$
acts transitively on it.  In summary,

\begin{lemma} \label{base_cycle} $K(z) = K_0(z)$; in particular it is 
a compact complex submanifold of $D$.
\end{lemma}

We write $C_0$ for the compact complex submanifold $K_0(z) \subset D$.
It will be the base cycle in a certain cycle space discussed below.  
The discussion leading
to Lemma \ref{base_cycle} shows that $C_0$ is both the unique $K$--orbit
in $D$ that is compact and the unique $K_0$--orbit in $D$ that is complex.
This is the origin of what is known as ``Matsuki duality.''

\begin{example}{\em
Let $Z$ be the complex projective space $\C\P^n$ and let $G_0
= SU(n,1)$.  Let $\{e_1 , \dots , e_{n+1} \}$ denote the standard basis 
of $\C^{n+1}$
relative to which the hermitian form defining $G_0$ is 
$\langle u , v \rangle$ = 
$\left ( \sum_{1 \leqq a \leqq n} u_a \overline{v_a} \right )$
$- u_{n+1} \overline{v_{n+1}}$\,.  Then $G_0$ has three orbits on $Z$: 
the (open) unit ball $\cB$ in $\C^n$ inside $Z$, consisting of the
negative definite lines, the $(2n-1)$--sphere $S$ which is the boundary
of $\cB$, consisting of the null lines, and the complement $D$ of $\cB \cup S$,
consisting of the positive definite lines.  
$D$ is the non--convex open $G_0$--orbit on $Z$.  Here $C_0$ is the hyperplane
at infinity, complement to $\C^n$ in $Z$.  In homogeneous coordinates
$[z^1, \dots , z^{n+1}]$, $\cB$ is given by 
$\sum_{1 \leqq a \leqq n}|z^a|^2 < |z^{n+1}|^2$, $S$ is given by 
$\sum_{1 \leqq a \leqq n}|z^a|^2 = |z^{n+1}|^2$, $D$ is given by
$\sum_{1 \leqq a \leqq n}|z^a|^2 > |z^{n+1}|^2$, and $C_0$ is given by
$|z^{n+1}|^2 = 0$.} \hfill $\diamondsuit$
\end{example}

Later we will see $C_0$ as the base cycle in $D$.  In this case 
$C_0$ is maximal among the (complex) subvarieties of $Z$ contained in $D$.

\section{Basic Facts on the Cycle Space.} \label{sec2}
\setcounter{equation}{0}

Basic facts about the cycle space are given in Wells \& Wolf \cite{WeW}
and in Wolf \cite{W7}.  We review 
some of that material now, and briefly indicate some of the applications of
cycle spaces to variation of Hodge structure, specifically to 
period matrix domains, and construction of automorphic cohomology classes
by Poincar\' e $\vartheta$--series.  Those applications, and several
others, will be discussed in more detail in Part III below using the 
tools which are described in Part II.
\medskip

{\bf Definition.}
Let $E = \{g \in G \mid gC_0 = C_0\}$.  Then $E$ is a closed complex subgroup
of $G$, so the quotient manifold
$\Omega := \{ gC_0 \mid g \in G \} \cong G/E$ has a natural 
structure of $G$--homogeneous complex manifold.  Since $C_0$ is compact and 
$D$ is open, the subset 
$\{ gC_0 \mid g \in G \text{ and } gC_0 \subset D \}$
is open in $\Omega$, and thus has a natural structure of complex 
manifold.  The {\em cycle space} of $D$ is
\begin{equation} \label{def_cyclespace}
\Omega_W(D): \text{ topological component of } C_0 \text{ in } 
\{ gC_0 \mid g \in G \text{ and } gC_0 \subset D \}.
\end{equation}
Thus $\Omega_W(D)$ has a natural structure of complex manifold.  
\medskip

\noindent {\bf Hermitian trichotomy.}
\medskip

In order to understand the structure of $Z$, $D$ and $\Omega_W(D)$ we may
assume that $G_0$ is simple, because $G_0$ is local direct product of simple
groups, and $Z$, $D$ and $\Omega_W(D)$ break up as global direct products
along the local direct product decomposition of $G_0$\,.  From this point
on $G_0$ is simple unless we say otherwise.

Since $\gg_0$ is simple and $\ge$ contains $\gk$, there are four 
possibilities, one trivial.  The trivial one is the case $\ge = \gg$,
in other words the case where $G_0$ acts transitively on $Z$, and $\Omega_W(D)$ 
is reduced to a single point. There are just a few possibilities for
this (Wolf \cite{W8}).  From now on we ignore this trivial case and concentrate 
on the other three:

\begin{enumerate}
\item {\sc Hermitian holomorphic case.} 
$G_0/K_0$ is a bounded symmetric domain $\cB$, we have the usual
$\gg = \gp_+ + \gk + \gp_-$\,, and $\ge$ is one of $\gk + \gp_\pm$\,. 
In this case $D$ is measurable, say $D = G_0/L_0$, and there is a
holomorphic double fibration\footnote{ 
By {\it holomorphic fibration} we mean a holomorphically locally trivial
fiber space, essentially a holomorphic fiber bundle except perhaps lacking
a complex structure group.}
\medskip

\setlength{\unitlength}{.08 cm}
\begin{picture}(180,18)
\put(80,15){$G_0/(L_0 \cap K_0)$}
\put(66,1){$D$}
\put(117,1){$\cB$}
\put(78,13){\vector(-1,-1){6}}
\put(107,13){\vector(1,-1){6}}
\end{picture}

\noindent
In other words, the two projections are simultaneously holomorphic for 
some choice between $\cB$ and the complex conjugate structure
$\overline{\cB}$ and some choice of invariant 
complex structure on $G_0/(L_0 \cap K_0)$.  In this case $G_0/(L_0 \cap K_0)$
is the incidence space $\cI(D)$ that we'll meet later, and $\Omega_W(D)$ is
$\cB$ or $\overline{\cB}$; see \cite{W7} or Wolf--Zierau \cite{WZ1}.   

\item {\sc Hermitian non--holomorphic case.} 
$G_0/K_0$ is a bounded symmetric domain $\cB$ and $\ge = \gk$. 
In this case we cannot adjust invariant complex structures so that the 
double fibration indicated just above will be holomorphic.  Here it was
recently proved that
$\Omega_W(D)$ is biholomorphic to $\cB \times \overline{\cB}$.  See
\cite{W7} or \cite{WZ1}, and Huckleberry--Wolf \cite{HuW3} or \cite{WZ3}.

\item  {\sc Generic (or non--hermitian) case.}
$G_0/K_0$ does not have a $G_0$--invariant complex structure.  Then
$\gk$ is a maximal subalgebra of $\gg$.  In this case $K$ is the identity
identity component of $E$, $\Omega$ is an affine homogeneous space, and
we will describe the structure of $\Omega_W(D)$ in that context.  The
precise structure was worked out only very recently in \cite{HuW3}
and Fels--Huckleberry \cite{FH}.

\end{enumerate}

\noindent
In the rest of this article, we describe complex geometric methods that
lead to the  developments indicated above, and to other applications and 
developments through the use of double fibration transforms.  The new
developments include aspects of the theory of holomorphic double fibration
transforms themselves (Section \ref{sec11}).  The areas 
of application include aspects of the representation theory of semisimple 
Lie groups (Section \ref{sec12}) and variation of Hodge structure 
(Section \ref{sec13}).
\medskip

\section{The Exhaustion Function and $(q+1)$--Completeness.} \label{sec3}
\setcounter{equation}{0}

Measurable open orbits $D = G_0(z) \subset Z$ carry an 
especially useful real analytic exhaustion 
function $\varphi : D \to \R$ whose Levi form $\cL(\varphi)$ has at least
$n-q$ positive eigenvalues at every point of $D$, where $n = \dim_\C D$ and
$q = \dim_\C C_0$\,.  Thus $\varphi$ is strongly $q$--pseudoconvex
and $D$ is $(q+1)$--complete.  In this section we review that development
from \cite{S1}, \cite{WeW} and \cite{SW}, and then we indicate applications
\cite{W7} to cohomology over $D$ and to the Stein property of $\Omega_W(D)$.
\medskip

The exhaustion function $\varphi : D \to \R$ was first described in
Schmid's thesis \cite{S1} in the setting where $G_0$ has a compact Cartan
subgroup and $Z = G/B$ where $B$ is a Borel subgroup of $G$.  The 
$\theta$--stable real form $G_u$ of $G$
acts transitively on $Z$.  The canonical line bundle $\K_Z \to Z$, and the
(dual) anticanonical line bundle $\K^*_Z \to Z$, are $G_u$--homogeneous
and have $G_u$--invariant metrics.  Let $h_u$ denote the $G_u$--invariant
hermitian metric on $\K^*_Z \to Z$.  In this setting the isotropy subgroup
$L_0$ of $G_0$ at a point $z \in D$ of the open orbit is just a compact Cartan
subgroup, so the anticanonical bundle $\K_D^* \to D$ has a $G_0$--invariant
hermitian metric $H_0$\,.  Then one has the $C^\omega$ (real analytic) 
positive function $\varphi = \log h_0/h_u$ on $D$.  If $g(z) \in \text{ bd}(D)$
then $\text{Ad}(g)(\gl + \gq_-) + \overline{\text{Ad}(g)(\gl + \gq_-)}
\subsetneqq \gg$, and it follows that $\varphi$ goes to infinity as one
approaches $g(z)$ from the interior of $D$.  From this one sees that
$\varphi$ is an exhaustion function for $D$.  Root space considerations
allow one to compute $\sqrt{-1}\partial \overline{\partial} \log h_0$ and
$\sqrt{-1}\partial \overline{\partial} \log h_u$ and see the Levi form
$\cL(\varphi)$ explicitly.  It follows immediately that $\cL(\varphi)$ has 
at least $n-q$ positive eigenvalues at every point of $D$.
\smallskip

Somewhat later, Wells and Wolf \cite{WeW} noted that Schmid's argument
could be adapted to the more general setting where the only requirement is
that the isotropy subgroup $L_0$ of $G_0$ at a point $z \in D$ is compact.
Somewhat after that, Schmid and Wolf \cite{SW} further adapted the argument
to the (even more general) situation where $D$ is a measurable open 
$G_0$--orbit in a complex flag manifold $Z = G/Q$.  Thus every measurable 
open orbit $D$ is $(q+1)$--complete.
\medskip

\noindent {\bf Cohomology vanishing theorems.}
\medskip

The theorem of Andreotti and Grauert
\cite{AnG} says that if a complex manifold $D$ is $(q+1)$--complete, and if
$\cS \to D$ is a coherent analytic sheaf, then the cohomologies 
$H^r(D;\cS) = 0$ for for all $r > q$.  
\smallskip

Since our measurable open $G_0$--orbit $D$ is $(q+1)$--complete, we have 
the vanishing $H^r(D;\cO(\E)) = 0$ for $r > q$, for every holomorphic vector
bundle $\E \to D$.  On the other hand, if
$\E \to D$ is a (sufficiently) negative bundle \cite{GrS} the methods
based on the Bott--Borel--Weil Theorem show that $H^r(D;\cO(\E)) = 0$ for 
$r < q$.  Thus, finally,
\begin{equation}
\text{if } \E \to D \text{ is (sufficiently) negative, then }
H^r(D;\cO(\E)) = 0 \text{ for } r \neq q.
\end{equation}
This will be very important when we discuss double fibration transforms.
\medskip

\noindent {\bf The Stein property for cycle spaces of measurable open orbits.}
\medskip

When $D$ is a measurable open orbit, Wolf \cite{W7} combined his extension of
boundary component theory of bounded symmetric domains (\cite{W2}, or see
\cite{W4}) with the exhaustion function $\varphi : D \to \R$, to prove that 
$\Omega_W(D)$ is a Stein manifold.  We review the argument.  
\smallskip

{\sc Case: $D$ is of hermitian holomorphic type.}  Here we may assume
$D = G_0(z)$ and $E = KP_-$, so $\Omega = G/E$ is the compact hermitian
symmetric space dual to the bounded symmetric domain $\cB$.  Thus
$\cB \subset \Omega_W(D) \subset \Omega$ and $\Omega_W(D)$ is invariant
by the action of $G_0$ on $\Omega$.  The $G_0$--orbit structure of $\Omega$,
and the closure relations among the $G_0$--orbits, are known precisely in 
terms of partial Cayley transforms (\cite{W2}, \cite{W4}).  If $\cO$
is a $G_0$--orbit in $\Omega_W(D)$ , then 
it contains every open $G_0$--orbit whose closure contains 
$\cO$, because $\Omega_W(D)$ is open in $\Omega$.  Some operator norm 
arguments show that $\Omega_W(D)$ cannot contain an open orbit different 
from $\cB$.  It follows that $\Omega_W(D) = \cB$, and in particular 
$\Omega_W(D)$ is Stein.
\smallskip

{\sc Case: $D$ is not of hermitian holomorphic type.}  Then $E$ has 
identity component $K$, so $E$ is reductive and $\Omega = G/E$ is affine.
Define $\beta: \Omega_W(D) \to \R^+$ by 
$\beta(gC_0) = \sup_{y \in C_0} \varphi(g(y))$.   Since $\varphi$ is an
exhaustion function and the $gC_0$ are compact, one sees that
$\beta : \Omega_W(D) \to  \R^+$ blows up at every boundary point of
$\Omega_W(D)$\,.  From the specific construction of $\varphi$, and a
close look at the real analytic variety given by $d\varphi = 0$, one sees 
that $\beta$ is continuous, piecewise 
$C^\omega$ and plurisubharmonic.  Now a modification
suggested by results of Docquier and Grauert \cite{DG} gives a $C^\omega$
strictly plurisubharmonic exhaustion function $\psi = \varphi + \nu$
constructed as follows.  Since $\Omega$ is Stein there is a proper
holomorphic embedding $f:\Omega \to \C^{2n+1}$ with closed image, by
Remmert's theorem.  Define $\nu(C) := || f(C)||^2$ for $C \in \Omega_W(D)$\,.
Since $\Omega_W(D)$ carries a strictly plurisubharmonic exhaustion function,
it is Stein.

\section{Early Problems and Results on the Double Fibration Transform.}
\label{sec4}
\setcounter{equation}{0}

We start this section with a review of some basic facts on the double 
fibration transform from \cite{WZ2}.  We then specialize (initially as in
\cite{WZ2}) to the case of an open orbit 
$D = G_0(z) \subset Z$, and present some new results that clear up several
open problems in that flag domain case.  Finally we give a quick indication of 
the consequences for
variation of Hodge structure and for semisimple representation theory.
Now we start with the general setup, indicate its technical
requirements, and specialize it to our flag domain situation.   Several of
the problems that come up here are settled later in Section \ref{sec11}
using methods developed in Part II below.
\medskip

\noindent {\bf Double fibration.}
\medskip

Let $D$ be a complex manifold (later it will be an open orbit of
a real reductive group $G_0$ on a complex flag manifold $Z = G/Q$
of its complexification ).  We suppose that $D$ fits into a {\sl 
holomorphic double fibration}, in other words that there are complex
manifolds $M$ and $\cI(D)$ with simultaneously holomorphic fibrations:
\begin{gather} \label{gen_doublefibration}
\setlength{\unitlength}{.08 cm}
\begin{picture}(180,18)
\put(81,15){$\cI(D)$}
\put(70,12){$\mu$}
\put(66,1){$D$}
\put(95,12){$\nu$}
\put(100,1){$M$}
\put(78,13){\vector(-1,-1){6}}
\put(92,13){\vector(1,-1){6}}
\end{picture}
\end{gather}
(Later $M$ will be a cycle space and $\cI(D)$ will be an incidence space for
points and cycles.)
Given a coherent analytic sheaf $\cE \to D$ we 
construct a coherent sheaf $\cE' \to M$ and a transform
\begin{equation}
P : H^q(D;\cE) \to H^0(M;\cE')
\end{equation}
under mild conditions on (\ref{gen_doublefibration}).  In fact we give 
several variations
on the construction.  This construction is fairly standard (see,
for example, \cite{BE}, \cite{PR1} and \cite{M}), but we need several
results specific to the case of flag domains.
\vfill\pagebreak

\noindent {\bf Pull--back.}  
\medskip

The first step is to pull cohomology back from $D$ to $\cI(D)$.  Let
$\mu^{-1}(\cE) \to \cI(D)$ denote the inverse image sheaf.  
For every integer $r \geqq 0$ there is a natural map
\begin{equation} \label{pull1}
\mu^{(r)} : H^r(D;\cE) \to H^r(\cI(D); \mu^{-1}(\cE))
\end{equation}
given on the \v Cech cocycle level by $\mu^{(r)}(c)(\sigma) = c(\mu(\sigma))$
where $c \in Z^r(D;\cE)$ and where
$\sigma = (w_0, \dots , w_r)$ is a simplex.  For $q \geqq 0$
we consider the Buchdahl $q$--condition
\begin{equation} \label{buchdahl_conditions}
\text{the fiber } F \text{ of } \mu: \cI(D) \to D \text{ is connected and }
H^r(F;\C) = 0 \text{ for } 1 \leqq r \leqq q-1.
\end{equation}

\begin{proposition} \label{buchdahl_theorem} {\rm (See \cite{Bu}.)}  
Fix $q \geqq 0$.
If {\em (\ref{buchdahl_conditions})} holds, then
{\rm (\ref{pull1})} is an isomorphism
for $r \leqq q-1$ and is injective for $r = q$.
If the fibers of $\mu$ are cohomologically acyclic then
{\rm (\ref{pull1})} is an isomorphism for all $r$.
\end{proposition}

As usual, $\cO_X \to X$ denotes the structure sheaf of a complex
manifold $X$ and $\cO(\E) \to X$ denotes the sheaf of germs of
holomorphic sections of a holomorphic vector bundle $\E \to X$.
Let $\mu^*(\cE) := \mu^{-1}(\cE) \widehat{\otimes}_{\mu^{-1}(\cO_D)} 
\cO_{\cI(D)} \to \cI(D)$ denote the pull--back sheaf.  It is a coherent analytic
sheaf of $\cO_{\cI(D)}$--modules.  If 
$\E = \cO(\E)$ for some holomorphic vector bundle
$\E \to D$, then $\mu^*(\cE) = \cO(\mu^*(\E))$,
where  $\mu^*(\E)$ is the pull--back bundle.  In any case,
$[\sigma] \mapsto [\sigma] \otimes 1$ defines a map
$i : \mu^{-1}(\cE) \to \mu^*(\cE)$
which in turn specifies maps in cohomology, the coefficient morphisms
\begin{equation}\label{pull2}
i_p : H^p(\cI(D); \mu^{-1}(\cE)) \to H^p(\cI(D); \mu^*(\cE)) \ \
\text{ for } p \geqq 0.
\end{equation}
Our natural pull--back maps are the compositions
$j^{(p)} = i_p \cdot \mu^{(p)}$ of (\ref{pull1}) and (\ref{pull2}):
\begin{equation} \label{pull}
j^{(p)} : H^p(D;\cE) \to H^p(\cI(D);\mu^*(\cE)) \ \ \text{ for } p \geqq 0.
\end{equation}
If $\cE = \cO(\E)$ for some holomorphic vector bundle
$\E \to D$, then $\mu^*(\cE) = \cO(\mu^*(\E))$, we
realize these sheaf cohomologies as Dolbeault cohomologies, and
the pull--back maps (\ref{pull}) are given by pulling back
$[\omega] \mapsto [\mu^*(\omega)]$ on the level of differential forms.
\medskip

\noindent {\bf Push--down.} 
\medskip

In order to push the $H^q(\cI(D);\mu^*(\cE))$ down to 
$M$ we assume that
\begin{equation} \label{proper_stein}
\nu : \cI(D) \to M \text{ is a proper map and }
M \text{ is a Stein manifold}.
\end{equation}
The Leray direct image sheaves
$\cR^p(\mu^*(\cE)) \to M$ are coherent 
\cite{GrR}.  As $M$ is Stein  
\begin{equation} 
H^q(M; \cR^p(\cE)) = 0 \ \ \text{ for } \ \ p \geqq 0 \text{ and } q > 0.
\end{equation}
Thus the Leray spectral sequence collapses and gives
\begin{equation} \label{collapse}
H^p(\cI(D); \mu^*(\cE)) \cong H^0(M;\cR^p(\mu^*(\cE))).
\end{equation}

\begin{definition} {\em
The {\em double fibration transform} for the holomorphic double 
fibration {\rm (\ref{gen_doublefibration})} is the composition 
\begin{equation} \label{transform}
P: H^p(D; \cE) \to H^0(M; \cR^p(\mu^*(\cE)))
\end{equation}
of the maps {\rm (\ref{pull})} and {\rm (\ref{collapse})}.
}
\end{definition}

In order that the double fibration transform (\ref{transform})
be useful, one wants two conditions to be satisfied.  They are
\begin{gather}
P: H^p(D; \cE) \to H^0(M; \cR^p(\mu^*(\cE))) \text{ should be injective, and}
\label{want_inj}\\
\text{there should be an explicit description of the image of } 
P \label{want_image}.
\end{gather}
Assuming (\ref{proper_stein}),
injectivity of $P$ is equivalent to injectivity of $j^{(p)}$ in
(\ref{pull}).  The most general way to approach this is the combination of
vanishing and negativity in Theorem \ref{general_injective} below, based on the
Buchdahl conditions (\ref{buchdahl_conditions}).  
\medskip

The general (assuming (\ref{proper_stein})) injectivity question uses a 
spectral sequence argument for the the relative Dolbeault complex
of the holomorphic fibration $\mu:\cI(D)\rightarrow D$.  See \cite{WZ2} for the
details.  The end result is

\begin{theorem} \label{general_injective} Fix $q \geqq 0$.  Suppose that 
the fiber $F$ of $\mu : \cI(D) \to D$ is connected satisfies 
{\rm (\ref{buchdahl_conditions})}.  Assume {\rm (\ref{proper_stein})} 
that $\nu : \cI(D) \to M$ is proper and $M$ is Stein, say with fiber $C$.  
Let $\Omega_\mu^r(\E) \to \cI(D)$ denote the sheaf of relative 
$\mu^*\E$--valued 
holomorphic $r$--forms on $\cI(D)$ with respect to $\mu : \cI(D) \to D$.
Suppose that $H^p(C; \Omega_\mu^r(\E)|_C) = 0$ for 
$p < q$, and $r \geqq 1$.  Then
$P :H^q(D;\cE) \to H^0(M; \cR^q(\mu^*\cE))$ is injective.
\end{theorem}

\begin{remark} {\em
In the cases of interest to us, $\cE=\cO(\E)$ for some
holomorphic vector bundle $\E \to D$, and $P$ has an
explicit formula.  The  Leray derived sheaf is given by
\begin{equation} \label{derived_bundle_gen}
\cR^q(\mu^*(\cO(\E))) = \cO(\E^\dagger) \text{ where } 
\E^\dagger \to M \text{ has fiber } 
H^q(\nu^{-1}(C); \cO(\mu^*(\E)|_{\nu^{-1}(C)}))
\text{ at } C.
\end{equation}
Let $\omega$ be an $\E$--valued
$(0,q)$--form on $D$  and $[\omega] \in H_{\overline{\partial}}^q(D,\E)$
its Dolbeault class.  Then 
\begin{equation*}
P([\omega]) \text{ is the section of }
\E^\dagger \to M \text{ whose value } P([\omega])(C) \text{ at }
C \in M \text{ is } [\mu^*(\omega)|_{\nu^{-1}(C)}].
\end{equation*}
In other words,
\begin{equation}
P([\omega])(C) = [\mu^*(\omega)|_{\nu^{-1}(C)}] \in 
H^0_{\overline{\partial}} (M; \E^\dagger). 
\end{equation}
This is most conveniently interpreted by viewing $P([\omega])(C)$ as
the Dolbeault class of $\omega|_{C}$, and by viewing $C \mapsto
[\omega|_{C}]$ as a holomorphic section of the holomorphic vector
bundle $\E^\dagger \to M$.
}
\end{remark}
\medskip

\noindent {\bf Flag domain case.}
\medskip

Now let $D = G_0(z_0)$ be an open orbit in the complex flag manifold 
$Z = G/Q$, and $M$ is replaced by the cycle space $\Omega_W(D)$.  Our double
fibration (\ref{gen_doublefibration}) is replaced by
\begin{gather} \label{flag_doublefibration}
\setlength{\unitlength}{.08 cm}
\begin{picture}(180,18)
\put(81,15){$\cI(D)$}
\put(70,12){$\mu$}
\put(66,1){$D$}
\put(95,12){$\nu$}
\put(100,1){$\Omega_W(D)$}
\put(78,13){\vector(-1,-1){6}}
\put(92,13){\vector(1,-1){6}}
\end{picture}
\end{gather}
where $\cI(D) := \{(z,C) \in D \times \Omega_W(D) \mid z \in C \}$ is the
incidence space.  Given a homogeneous holomorphic vector bundle
$\E \to D$, and the number $q = \dim_\C C_0$\,, the Leray derived sheaf
involved in the double fibration transform satisfies (\ref{derived_bundle_gen}).
Here that takes the form
\begin{equation}\label{derived_bundle}
\cR^q(\mu^*(\cO(\E))) = \cO(\E^\dagger) \text{ where } 
\E^\dagger \to \Omega_W(D) \text{ has fiber } H^q(C; \cO(\E|_{C}))
\text{ at } C \in \Omega_W(D).
\end{equation}
Evidently, $\E^\dagger \to \Omega_W(D)$ is globally $G_0$--homogeneous and
infinitesimally $\gg$--homogeneous, and $H^q(C; \cO(\E|_{C}))$ can be 
calculated in any reasonable case from the Bott--Borel--Weil Theorem,
especially when $\E \to D$ is negative.  Thus $\cR^q(\mu^*(\cO(\E)))$ is
given explicitly by (\ref{derived_bundle}) in the flag domain case.
\medskip

Using methods of complex geometry to be described in Part II, we will
see in Part III that $\Omega_W(D)$ is a contractible Stein
manifold in general, so $\E^\dagger \to \Omega_W(D)$ is holomorphically
trivial, and that $F$ satisfies (\ref{buchdahl_conditions}) for all $q$,
so the double fibration transforms are injective.  Thus, in the flag
domain case, we will have a complete answer to (\ref{want_inj}) and some
sharp progress toward (\ref{want_image}).  See Section \ref{sec11}.
\medskip

Now let us take a quick historical look back at the hermitian trichotomy
of Section \ref{sec2}.  In the hermitian holomorphic case 
$\Omega_W(D)$ is $\cB$ or $\overline{\cB}$, and one knows \cite[Section 4]{WZ2}
that $F$ and $\Omega_W(D)$ are contractible Stein manifolds.  In the 
hermitian nonholomorphic case, where $\Omega_W(D) = \cB \times \overline{\cB}$
(\cite{W7} or \cite{WZ1}, and \cite{HuW3} or \cite{WZ3}), there had only
been partial information (see \cite[Theorem 6.6]{WZ2}) on contractibility 
of $F$.  There had been essentially no information in the nonhermitian case.
\smallskip

One more remark.  In some cases one knows that $H^q(D; \cE)$
is an irreducible representation space for a group under which all
our constructions are equivariant, and one sees directly that $P$ is an 
intertwining operator, thus zero or injective.  In practice, however,
we usually look for implications in the other directions.  See Section 
\ref{sec11}.
\medskip 

\bigskip
\centerline{\Large \bf \underline{Part II: The Complex Geometric Approach.}}
\bigskip

In this Part we describe the methods and results in complex geometry that
lead to a structure theory for the cycle space, various associated
domains, and the double fibration transform.

\section{Introduction to the Complex Geometric Approach.} \label{sec5}
\setcounter{equation}{0}
Our goal here is to explain recent results which have led to the
characterization of $\Omega _W(D)$ as being equivalent to a 
certain universal domain $\Omega _{AG}$ in all but the well--understood
trivial and the hermitian holomorphic cases as discussed in Section \ref{sec2}.
Without further reference we exclude those cases in the sequel.

\medskip
In the present section we outline the relevant results and methods
in a nontechnical way.  In the following sections we give enough details 
so that the reader should have no difficulty working through
the literature.  See \cite {HuW3} and \cite {FH} for complete details.

\medskip
Building on experience with transversal varieties gained 
in \cite {HuW1}, \cite {W9} and \cite {HuS}, the Schubert domain
$\Omega _S(D)$ was introduced in \cite {H} as a tool for 
understanding complex analytic properties of $\Omega _W(D)$. The
motivation for this is quite transparent in \cite {HuS}, although
there only the case of $G_0=SL_n(\mathbb R)$ was considered.  
\smallskip

We discuss \cite {HuS} with the
benefit of hindsight and a more up to date notation.
Let $G_0=K_0A_0N_0$ be an Iwasawa decomposition of the real form
$G_0$. The corresponding set $KAN$ corresponds to the open $AN$--orbit
in the spherical affine homogeneous space $\Omega :=G/K$ and is
a proper, Zariski open subset of $G$. We refer to a Borel subgroup
$B \subset G$ as an {\it Iwasawa--Borel subgroup} of $G$ if it contains
an Iwasawa factor $A_0N_0$.  Of course these are just
the Borel subgroups which occur as the isotropy groups at points
of the closed $G_0$--orbit in $G/B$. Given an Iwasawa--Borel subgroup 
$B$ in $G$ and a point $z\in Z=G/Q$, the closure $S=\text{c}\ell (\cO)$ 
of the orbit $\cO=B.z$ is the associated {\it Iwasawa--Schubert variety}
or {\it Schubert cycle}. Let
$Y$ denote the complement of the {\it Schubert cell} $\cO \in S$, i.e.,
$Y:=S\setminus \cO$
\smallskip.

An Iwasawa--Schubert variety $T = \text{c}\ell (\cO)$ is {\it transversal} 
(relative to an open $G_0$--orbit $D$) if
(1) $T\cap D$ is nonempty and contained in the open $B$--orbit $\cO$, and
(2) codim\,$T=q= \dim\,C_0$ and the intersection $T\cap C_0$ is
transversal at each of its points.
In this dimension $n-q$, the Iwasawa decomposition
of $G_0$ implies that $Y\subset Z\setminus D$ (Theorem \ref{codimension}).
\smallskip

It was shown in \cite {HuS} that for an open $SL_n(\mathbb R)$--orbit 
$D$ in a flag manifold $Z$ of $SL_n(\mathbb C)$, and 
$C\in \text{bd}(\Omega _W(D))$, there exists a transversal
Schubert variety $T$ such that $Y\cap C\not =\emptyset$.  The
method of incidence varieties (\cite {BK},\cite {BM}) then shows
that the algebraic variety $A_Y:=\{C\in \Omega :C\cap Y\}$
contains the polar set $H_Y$ of a nonconstant meromorphic function
which is produced by the method of trace transform. In this context
we refer to $H_Y$ as an {\it incidence hypersurface}. 
It follows immediately that $\Omega _W(D)$ is a Stein
domain in $\Omega $, because each of its boundary points is contained
in an analytic hypersurface which is entirely contained in its complement.
\medskip

Returning to the general case, we note that if $T$ is Poincar\'e dual
to $C_0$, then it is indeed a transversal Schubert variety. However,
at least initially there is no reason to believe that, given
$C\in \text{bd}(\Omega _W(D))$, there exists 
$T$ with $C\cap Y\not =\emptyset $.
\medskip

Domains with supporting analytic hypersurfaces at each
of their boundary points have optimal character from the complex
analytic viewpoint. Thus the following Schubert domain was introduced 
in \cite {H}. Let $\Omega _S(D)$ be the connected component
containing $C_0$ of the complement of the union of all Iwasawa--Borel invariant
intersection hypersurfaces $H_Y$ which are defined by transversal
Schubert varieties.  Since $K_0$ acts transitively on such Borel groups,
the set which is removed is a compact family of hypersurfaces and
clearly $\Omega _S(D)$ is a proper, open, Stein domain in $\Omega $
which contains $\Omega _W(D)$.  It should be emphasized that
the sets $Y$ are possibly very far away from the boundary
of the cycle space and that the inclusion 
$\Omega _W(D)\subset \Omega _S(D)$ could theoretically be proper.
\smallskip

In principle there could be a plethora of domains $\Omega _S(D)$, but
experience with real forms of $SL_n(\mathbb C)$
(see \cite {HuW2} for the remaining cases) and classical hermitian
symmetric spaces (\cite {N1}, \cite {WZ1}) suggests that 
they might all be the same, agreeing with a domain which is
defined by removing {\it all} Iwasawa--Borel invariant hypersurfaces
from $\Omega $.
\smallskip

More precisely, let $\Omega _I$ be defined as the connected component
containing $C_0$ of the complement of the union of all $B$--invariant
algebraic hypersurfaces in $\Omega $, where $B$ runs over all 
Iwasawa--Borel subgroups of $G$.  Note that, without further information,
$\Omega _I$ could theoretically be empty. But in any case, using the
same argument as above, it is a Stein domain and of course
$\Omega _I\subset \Omega _S(D)$ for every open orbit $D$ in
every $G$--flag manifold $Z=G/Q$.
In order to understand the relation of the cycle spaces 
to the universal domain $\Omega _{AG}$, it is natural
to compare $\Omega _I$ and $\Omega _{AG}$.  In fact, using the
identification of $\Omega _{AG}$ with the maximal domain 
of definition $\Omega _{adpt}$ of the {\it adapted complex
structure} in the tangent bundle of the Riemannian symmetric
space $G_0/K_0$ (\cite{BHH},\cite {Ha}; see Section \ref{sec6} below),
an elementary argument involving the transported norm function
shows that $\Omega _{AG}\subset \Omega _I$ \cite {H}.
\smallskip

In another guise $\Omega _I$ had been considered and, 
from a completely different viewpoint the above
inclusion had been shown for classical groups \cite {KS}.
Recently, a purely algebraic proof was given in \cite {M2}.
\smallskip

Regarding $\Omega _I$ as the {\it polar} $\widehat X_0$ 
(the two definitions are easily seen to be equivalent (see
\cite{H}), but nevertheless reflect two very different aspects 
of the subject), Barchini proved the opposite inclusion \cite {B}.
\smallskip

As a consequence of considerations of cycle spaces
associated to non--open $G_0$--orbits, it is implicitly shown
by case by case methods in \cite {GM} that 
$\Omega _{AG}\subset \Omega _W(D)$ for
classical groups and exceptional hermitian groups.  Thus
$\Omega _{AG}=\Omega _I\subset \Omega _S(D)$
was known at this point.
\medskip

We indicate our contributions; they will be sketched in more
detail below.  In \cite {HuW3} we carried out the general program
which was indicated by the naive flag arguments of \cite {HuS}:
For every $C\in \text{bd}(\Omega _W(D))$ there exists a transversal
Schubert variety $T$ so that $Y\cap C\not= \emptyset$; in particular,
$\Omega _W(D)=\Omega _S(D)$ and thus $\Omega _W(D)$ is
a Stein domain with respect to functions in the image of the
trace transform.
\smallskip

The construction of $T$ results from the construction of
a Schubert variety $Y_p$ containing an arbitrarily given boundary point 
$p\in \text{bd}(D)$, with codim\,$Y_p=q+1$, and of course
$Y_p\subset Z\setminus D$.  This is a consequence of a sort of 
triality that initiates with Matsuki duality.  It is of interest
that the moment map, Morse--theoretic method for realizing this duality
(\cite {MUV},\cite {BL}) plays a direct role in our considerations
(see Section \ref{sec7}). Furthermore, it is shown that, given $C$ in the boundary
of the cycle domain, there exits $p\in C$ so that $Y_p$ is
indeed the $Y$ of a transversal Schubert variety $T$ (see Theorem \ref{trans}).
So the incidence variety method along with the inclusion
$\Omega _{AG}\subset \Omega _I$ yields 
$\Omega _{AG}\subset \Omega _W(D)= \Omega _S(D)$.
\medskip

Satisfied with identifying $\Omega _{AG}$
with $\Omega _{adpt}$ or $\Omega _I$,
we have actually not given its original definition \cite {AkG}.
We do that now because we will need its computable nature.
\medskip

The symmetric space $\Omega _0:=G_0/K_0$ is embedded
as a totally real submanifold of half--dimension of $\Omega =G/K$ as
the $G_0$--orbit of the neutral point. General principles
imply that there are $G_0$--invariant neighborhoods of $\Omega _0$
in $\Omega $ on which the $G_0$--action is proper.  
One looks for a canonically defined domain for which this is the case.
\smallskip

Let ${\mathfrak g}_0={\mathfrak k}_0\oplus {\mathfrak p}_0$
be a Cartan decomposition set up in the usual way with respect to
a compact real form ${\mathfrak g_u}$ of ${\mathfrak g}$.  Let 
${\mathfrak a}_0\subset {\mathfrak p}_0$ be maximal abelian as above.  For
$\alpha $ a  root of ${\mathfrak a}_0$, let 
$H_\alpha :=\{A\in {\mathfrak a}_0: \alpha =\frac {\pi}{2}\}$
and define
$\omega _0$ as the connected component containing $0$ of the
complement in ${\mathfrak a}_0$ of the union of the $H_\alpha $.
Then $\Omega _{AG}:=G_0.\exp(i\omega _0)z_0$, where $z_0$ is
a base point corresponding to $C_0$, is an open neighborhood
of $\Omega _0$ and is maximal with respect to the property that
every $G_0$--isotropy group is compact.  The $G_0$--action on
$\Omega _{AG}$ is proper (see Section \ref{sec6}).
\medskip

The fact that the $G_0$--action is proper indicates that
an invariant metric (perhaps of canonical nature) is playing
a role. Thus, in \cite {H} methods were introduced to study
the {\it hyperbolicity} of such domains.  Let us recall the
basic facts which are relevant for such considerations.
\medskip

The Kobayashi pseudo--metric on a complex manifold $X$ can be
defined as follows. First, define a disk in $X$ as the
biholomorphic image of the unit disk in the complex plane.
A chain of disks is the union of finitely many such disks
which overlap (on open subsets) to form a connected set.
Given $p,q\in X$ consider a chain 
$\kappa =\Delta _1\cup \ldots \cup \Delta _m$ with $p\in \Delta _1$,
$q\in \Delta _m$ and $\Delta _i\cap \Delta _{i+1}\not =\emptyset $.
Let $p_i\in \Delta _i\cap \Delta _{i+1}$ and $d_i$ be the distance from
$p_{i-1}$ to $p_i$ computed in the Poincar\'e metric of 
$\Delta _i$.  Adding up these distances we obtain a number
$d(\kappa )$ which also depends on the choices of the $p_i$ which
are regarded as part of the data of the chain.  
Finally, define the Kobayashi pseudo--distance
between $p$ and $q$ as the minimum of all such $d(\kappa )$ as
$\kappa $ runs over all such chains.  This defines a pseudo--distance
function on $X\times X$ which, if it is nonzero for all
$p,q\in X$, is the {\it Kobayashi metric} \cite {K}. In this
case $X$ is called {\it Kobayashi hyperbolic}.

\medskip
The methods introduced in  \cite {H} could be regarded as leading to the
group--theoretic version of the fact that the complement in
$\mathbb C\mathbb P_m$ of $2m+1$ hyperplanes in general position
is Kobayashi hyperbolic (a result of classical geometry).  It is
shown in \cite {H} that, e.g., $\Omega _I$ is hyperbolic.
These methods, refined in \cite {FH}, lead to the
following result (see Theorem \ref{k_hyper}) in the non--hermitian case.
Let $H$ be any Iwasawa--Borel invariant hypersurface in $\Omega $
and define $\Omega _H$ to be the connected component containing
$C_0$ in $\Omega $ of the complement of $\bigcup_{k\in K_0}k(H)$.
In other words, the definition of $\Omega _H$ is analogous to $\Omega _S(D)$
and $\Omega _I$ except that one initially has only
a single hypersurface.  

\begin {theorem}
If $G_0$ is not of hermitian type, then $\Omega _H$ is Kobayashi hyperbolic.
\end {theorem}

Summarizing the above, if $G_0$ is not of hermitian type then, for any 
Iwasawa--Borel invariant hypersurface $H$ in the complement of $\Omega _W(D)$
\begin {equation}\label {inclusions}
\Omega _{AG}=\Omega _I\subset \Omega _W(D)=\Omega_S(D)\subset \Omega _H
\end {equation}
and each of these domains is Stein and Kobayashi hyperbolic.
\medskip

The main new development in \cite {FH} is summarized as follows.

\begin {theorem} \label {FH-main}
Suppose that $G_0$ is not of hermitian type.
If $\widehat \Omega $ is a $G_0$--invariant, Kobayashi hyperbolic
Stein domain that contains $\Omega _{AG}$, then 
$\widehat \Omega =\Omega _{AG}.$
\end {theorem}
It follows that the inclusions in (\ref {inclusions}) are equalities.
With some additional remarks to handle the hermitian
case (which in fact is much simpler), it follows that
$\Omega _W(D)=\Omega _{AG}$ in all but the well understood trivial and 
hermitian holomorphic cases discussed in Section \ref{sec2}.
\medskip

We complete these introductory remarks by outlining the
basic ideas of \cite {FH}. Roughly speaking, the goal of that work
is to reduce to the case of $SL_2$.  This requires a 
rather detailed analysis of the $G_0$--action on
$\text{bd}(\Omega _{AG})$, in particular
a detailed description of closed orbits and orbit closures
(see the subsection on closed orbits in Section \ref{sec10} below).  
The boundary $\text{bd}(\Omega _{AG})$ 
is not smooth, but for the purposes of \cite {FH}
it is sufficient to consider points in a tractable open, dense stratum
$\text{bd}_{gen}(\Omega _{AG})$ (see the subsection on genericity in
Section \ref{sec10} below).  For 
$z\in \text{bd}_{gen}(\Omega _{AG})$ one determines a $3$--dimensional
simple subgroup $S \subset G$ which is defined over $\mathbb R$ such
that $S.z\cap \Omega _{AG}$ is the associated domain $\Omega _{AG}(S)$
for the group $S$ with respect to its noncompact real form $S_0$.
\smallskip

Since no difficulties are introduced by going to finite covers,
it may be assumed that the associated affine variety $\Omega (S)$
is the complement of the diagonal in 
$\mathbb C\mathbb P_1\times \mathbb C\mathbb P_1$, i.e., the 
$2$--dimensional affine quadric. Regarding it as a closed $S$--orbit in
$\Omega $, we refer to $\Omega (S)$ as a {\it $Q_2$--slice}.
\smallskip

The $S$--action on $\mathbb C\mathbb P_1\times \mathbb C\mathbb P_1$
is the standard diagonal action.  The intersection
$\Omega (S)\cap \Omega _{AG}=\Omega _{AG}(S)$ is a $2$--dimensional
polydisk $\Delta $ which can be regarded as
being the product $\cB^+\times \cB^-$ of an $S_0$--invariant ($1$--dimensional)
disk with its exterior in $\mathbb C\mathbb P_1$. Due to the 
genericity assumption, the point $z$ can be chosen to lie
in $\cB^+\times \text{bd}(\cB^-)$.
\smallskip

It is shown in \cite{FH} that if $\widehat \Omega(S)$ is a $S_0$--invariant 
Stein domain in $\Omega (S)$ which contains $\Omega _{AG}(S)$ and
the boundary point $z$, then it contains the open set
$\Omega \cap (\cB^+\times \mathbb C\mathbb P_1) \cong B^+\times \mathbb C$.
It is straightforward to check that if a complex manifold $X$ contains
a (biholomorphic) copy of $Y$ of $\mathbb C$, then the Kobayashi
distance between any two points in $Y$ is zero.  Thus the domain
$\widehat \Omega (S)$ is certainly not Kobayashi hyperbolic.
\smallskip

The proof of Theorem \ref {FH-main} is then immediate, because, if
$\widehat \Omega$ is a $G_0$--invariant, Stein domain which 
properly contains $\Omega _{AG}$, then there is a generic
boundary point $z$ such that $\widehat \Omega (S)\cap \widehat \Omega $
contains both $\Omega _{AG}(S)$ and $z$.  Since $\Omega _{AG}(S)$
is Stein and $S_0$--invariant, it contains the copies of $\mathbb C$
as above and therefore $\widehat \Omega $ contains these as well.
As a consequence $\widehat \Omega $ is not hyperbolic and
Theorem \ref {FH-main} follows.
  
\section{The Equivalences $\Omega _{adpt}=\Omega _{AG}=\Omega _I$.} \label{sec6}
\setcounter{equation}{0}

Here three $G_0$ domains are introduced from three different viewpoints.
The domain $\Omega _{adpt}$ in the tangent bundle of the Riemannian
symmetric space $\Omega _0=G_0/K_0$ can be defined by either metric
or symplectic properties.  The equivalence of these two ways of
viewing $\Omega _{adpt}$ are important for complex analytic considerations.
\medskip

The domains $\Omega _{AG}$ and $\Omega _I$ are defined as neighborhoods
of $\Omega _0$ in the affine homogeneous space $\Omega =G/K$, where
a base point $x_0$ has been chosen so that $\Omega _0=G_0.x_0$.  The
domain $\Omega _{AG}$ is defined from the point of view of group actions. 
The domain $\Omega _I$ can be seen from several viewpoints.  Ours
is that of incidence divisors which are defined by
Schubert varieties in $Z$ of Iwasawa--Borel subgroups $B \subset G$.
\medskip

\noindent {\bf Adapted complex structures.} 
\medskip

Beginning with the metric standpoint, let $(M,g)$ be a real analytic
Riemannian manifold which for simplicity is assumed to be complete.
The differential $\gamma _*$ of a geodesic is a map
$\gamma _*:T\R\to TM$ of tangent bundles and can be viewed as an orbit of the
$\R$--action defined by geodesic flow.  Let $\R^*$
act by scalar multiplication in the fibers of $TM$.  This action
extends in the usual way to a map $\R\times TM\to TM$.
Identify $T\R$ with the complex plane $\C$ by
$(t,s{\frac {d}{dt}}|_t)\mapsto t+is$ and define
$\gamma ^\C:\C\cong T\R\to TM$
by $z=s+it\mapsto s\cdot \gamma _*(t)$, where the multiplication
comes from the $\R^*$--action.

\begin{definition}{\em
An integrable complex structure $J$ on a starlike neighborhood
of $\cA$ of the $0$--section in $TM$ is 
{\it adapted} if for every geodesic $\gamma $ there is a
disk $\Delta =\Delta (0)\subset \C$ such that
$\gamma ^\C|_\Delta :\Delta \to \cA$ is
holomorphic.
}
\end{definition}
The existence and uniqueness of adapted structures are proved
in (\cite {LS}; also see \cite {Ha}).  The uniqueness statement
says: 
If $J_1$ and $J_2$ are adapted structures on $\cA$, then $J_1=J_2$.

The symplectic side of the picture was developed in
(\cite {GuS}; also see \cite {Bu}).  Let
$\lambda _{std}$ (resp. $\omega _{std}$) be the standard
$1$--form (resp. $2$--form) on the cotangent bundle
$T^*M$.  Let $\rho_g =TM\to \mathbb R$ be the norm--function
defined by the metric, i.e., $\rho _g(\cdot ) :=\Vert \cdot \Vert^2_g$.
If $\psi _g:TM\to T^*M$ denotes the diffeomorphism defined by the
metric, then we have the forms $\theta _g:=\psi ^*_g(\lambda _{std})$
and $\omega _g:=\psi^*_g(\omega _{std})$.
Given ${\cal A}$ as above, $d^c\rho _g=\theta _g$ is regarded as
a differential equation for a complex structure $J$ on ${\cal A}$.
The local existence of integrable such structures is shown in
\cite {GuS} and the same strong uniqueness theorem as that stated above
is proved; in particular, the locally defined $J$'s automatically glue
together. Furthermore,
$dd^c\rho _g=d\theta _g=\omega _g$
is K\"ahlerian, i.e., $\rho _g$ is a strictly plurisubharmonic
function on every adapted neighborhood ${\cal A}$.

\smallskip
The connection between these two notions of adapted structure
is given by the following result.

\begin {theorem} {\rm \cite {LS}}
The $1$--form $\theta _g$ on a domain ${\cal A}$ equipped with
the adapted complex structure in the Riemannian sense satisfies
$d^c\rho _g=\theta _g$.
\end {theorem}

As a consequence of the uniqueness theorem, the
Riemannian and symplectic notions of adapted structure are
equivalent.  This allows the use of properties of
plurisubharmonic functions, which we now briefly summarize,
in the Riemannian setting.
\medskip

\noindent {\bf Basic properties of plurisubharmonic functions.}
\medskip

A (smooth) strictly plurisubharmonic function $\rho :X\to \mathbb R$ on
a complex manifold is by definition a potential of a
K\"ahler form $dd^c\rho =\omega $. In other words, in holomorphic
coordinates the complex Hessian
$H(\rho ):=
\begin {pmatrix}
\frac{\partial ^2\rho}{\partial z_i\partial \bar z_j}
\end {pmatrix}$
is positive definite.  Equivalently, the restriction $\rho|_C$
to every (local) complex curve $C$ is strictly subharmonic in the
sense that the Laplacian of $\rho|_C$ is negative.  On disks
in $C$ such functions have the (strong) mean value property.
\smallskip

Strictly plurisubharmonic functions have strong convexity properties.
In fact, holomorphic coordinates can be chosen so that $\rho $
is strictly convex on the underlying real domain.  Thus the
maximum principle holds: A strictly plurisubharmonic function never
takes on a (local) maximum value.  The following also reflects
this strong convexity.

\begin {proposition} \label {psh}
Let $\rho :X\to \mathbb R$ be strictly plurisubharmonic.  Let
$M\subset X$ be a connected local real submanifold such that
{\rm (1)} $\rho$ has some constant value $c$ on $M$ and
{\rm (2)} $c$ is a minimal value of $\rho $ in a neighborhood of $M$.
Then $\dim_\mathbb RM \le \dim_\mathbb CX$
\end {proposition}
This follows immediately from the fact that the complex Hessian of $\rho$
is positive definite, and therefore if $V$ is a real subspace of the
tangent space that is isotropic with respect to the real Hessian of $\rho$\,,
then $\dim_\R V \leqq \dim_\C X$.
\medskip

\noindent {\bf The adapted structure for Riemannian symmetric spaces.}
\medskip

Here only irreducible Riemannian symmetric spaces $\Omega _0=G_0/K_0$
of negative curvature are considered.  In other words $G_0$ is
noncompact and simple.  If ${\mathfrak g}_0=
{\mathfrak k}_0\oplus {\mathfrak p}_0$ is a Cartan decomposition,
then $T\Omega _0=G_0\times _{K_0}{\mathfrak p}_0$.  
Now we recall the Riemannian
notion of adapted complex structure on domains ${\cal A}$ in
$T\Omega _0$.
\medskip

Define the {\it polar coordinate map} $\Psi:T\Omega_0 \to \Omega =G/K$ by
$[(g_0,\xi)]\mapsto g_0\exp(i\xi).x_0$. 
It is well defined and $G_0$--equivariant.  Furthermore,
$\Psi_*(p):T_pT\Omega _0\to T_{\Psi(p)}\Omega $ is an
isomorphism at points $p$ of the $0$--section.

\medskip
Let $\Omega _{adpt}$ be the connected component containing
the $0$--section of the set where $\Psi$ has maximal rank, which is
$\{ p\in T\Omega _0 \mid {\Psi}_*(p) \text{ is an isomorphism}\} $.
Let $J$ be the (integrable) complex structure on
$\Omega _{adpt}$ which is defined by pulling back the complex
structure of $\Omega $ by ${\Psi}$.

\begin {proposition}
The structure $J$ on $\Omega _{adpt}$ is adapted.
\end {proposition}

\noindent {\bf Proof.}
Let $x_0 = (x_0,0)$, the neutral point in $T\Omega_0$\,.  Identify $\Omega _0$
with $G_0.x_0$. Recall that the geodesics through $g(x_0)$ are given by
$1$--parameter groups: $\gamma (t)=g \exp(t \xi).x_0$ for 
$\xi\in {\mathfrak p}_0$\,.  Thus
$s\gamma _*(t)=[(g.\exp(t\xi).e, s\xi)]$
and ${\Psi}\circ \gamma ^\mathbb C(t+is)=g\exp((t+is)\xi )x_0$.
Consequently, for an appropriately small disk,
$\gamma ^\mathbb C\vert \Delta:\Delta\to \Omega _{adpt}$ is
holomorphic. \hfill $\square$
\medskip

\noindent {\bf Proper actions.}
\medskip

An action $L\times M\to M$ of a topological group on a topological
space is said to be {\it proper} if the induced map
$L\times M\to M\times M$, $(g,x)\mapsto (g(x),x)$, is proper.
Under minimal assumptions on the spaces at hand, this can be
expressed as follows: For all sequences $\{g_n\}\subset L$
and $\{ x_n\}\subset M$ such that $x_n\mapsto x$ and
$g_n(x_n)\mapsto y$ there exists a convergent subsequence
$g_{n_k}\mapsto g$ in $L$. All isotropy groups of a proper
action are compact.

\medskip
The $G_0$--action on $\Omega _0=G_0/K_0$ is proper, so
the $G_0$--action on $M=T\Omega _0$ is as well.  Although
${\Psi}:\Omega _{adpt}\to \Omega $ has
finite fibers, it does not immediately follow that the
$G_0$--action on its image is proper.  Nevertheless, the
$G_0$--isotropy groups in ${\Psi}(\Omega _{adpt})$ are at most
finite extensions of the corresponding (compact) isotropy groups
in $\Omega _{adpt}$ and therefore are themselves compact.
It would therefore be natural to consider canonically defined neighborhoods
of $\Omega _0$ in $\Omega $ in which the $G_0$--isotropy groups
are compact.

\medskip
From the point of view of the Riemannian conjugate locus there is
a very natural candidate for such a domain.  See \cite {C}.  In
order to define it, assume as usual that the Cartan involution
of $G_0$ is the restriction of that for $G$ which in turn
defines its maximal compact subgroup $G_u$.  Consider the
restriction of the polar coordinate map ${\Psi}$ to
the fiber $T_{x_0}\Omega _0\cong{\mathfrak p}_0$ at the neutral point
in $T\Omega _0$. It maps sufficiently small open neighborhoods
of $0\in {\mathfrak p}_0$ diffeomorphically onto neighborhoods
of the neutral point $x_0$ in the $G_u$--orbit $G_u.x_0=G_u/K_0$
in $\Omega $.

\medskip
The maximal such set is determined as follows by the
conjugate locus of the invariant metric.  Let ${\mathfrak a}_0$ be a
maximal abelian subalgebra of ${\mathfrak g}_0$ which
is contained in ${\mathfrak p}_0$, $\Lambda _0$ be its
set of (real) roots, and for $\alpha \in \Lambda _0$, let
$H_\alpha $ be the affine hypersurface
$\{ \xi \in {\mathfrak a}_0: \alpha (\xi )=\frac {\pi}{2}\} $.
Define $\omega _0$ to be the connected component
containing $0\in {\mathfrak p}_0$ of
${\mathfrak a}_0\setminus  \bigcup_{\alpha \in \Lambda _0}\ H_\alpha $
and let $\Sigma _0:=K_0 \exp(\omega _0)$.  Then
$G_0\times _{K_0}\Sigma _0$ is naturally embedded
in $T\Omega _0$ as an open neighborhood of the $0$--section by the
action map $(g,\xi)\mapsto [(g_0,\xi)]$ and
${\Psi}|_{(G_0\times _{K_0}\Sigma _0)}:G_0\times _{K_0}\Sigma \to \Omega$
is a diffeomorphism onto its image (see e.g. \cite {C}).
This image ${\Psi}(G_0\times _{K_0}\Sigma _0)=G_0.\exp(i\omega _0).x_0$
was considered in \cite {AkG}, and we denote
by $\Omega _{AG}$\,.  The situation can now be summarized as follows.

\begin {proposition} \label {easy inclusion}
The restriction of ${\Psi}$ to the open subset $G_0\times _{K_0}\Sigma _0$
of $\Omega _{adpt}$ is a diffeomorphism onto its image $\Omega _{AG}$; in
particular the $G_0$--action on $\Omega _{AG}$ is proper.  Furthermore,
$\Omega _{AG}$ comes equipped with the ${\Psi}$--induced, $G_0$--invariant
strictly plurisubharmonic function $\rho =\rho _g\circ {\Psi}^{-1}$
and its associated K\"ahler form $\omega :=dd^c\rho $.
\end {proposition}

The following result sheds more light on the picture.

\begin {theorem} {\rm \cite {BHH}}
The domain $G_0\times _{K_0}\Sigma _0$ is a maximal domain of
definition for the adapted complex structure.  In particular
$\Omega _{adpt}=G_0\times _{K_0}\Sigma _0$ and
${\Psi}:\Omega _{adpt}\to \Omega _{AG}$ is biholomorphic.
\end {theorem}
\medskip

\noindent {\bf Incidence geometry and the domain $\Omega _I$\,.}
\medskip

Recall that an Iwasawa--Borel subgroup is a
Borel subgroup $B \subset G$ that contains a factor $A_0N_0$ of an Iwasawa 
decomposition $G_0=K_0A_0N_0$\,, where $K_0$ can be any maximal
compact subgroup of $G_0$.

\medskip
Consider our usual cycle space setup for an open $G_0$--orbit
$D$ in $Z=G/Q$, where $\Omega _W(D)$ is regarded as an open
neighborhood of $\Omega _0=G_0.C_0=G_0.x_0$ in $\Omega =G/K$.
Recall that a $B$--Schubert variety $S$ in $Z$ is defined
to be the closure of a $B$-orbit $\cO$ in $Z$.  Write
$S=\cO \cup Y$, where $Y$ is the (finite) union of $B$--orbits
on the boundary of $\cO$. 
The following remark motivates a number of our considerations.

\begin {theorem}  \label {codimension}
If $S\cap D\not =\emptyset $, then
{\rm codim}$_Z S \leqq \dim_\mathbb C C_0$.
\end {theorem}

\noindent {\bf Proof.}
If $S\cap D \not =\emptyset $, then $\cO\cap D\not=\emptyset $ as
well.  Without loss of generality we may assume that $G_0=K_0A_0N_0$,
$B \supset A_0N_0$ and $K_0.z_0=C_0$ is the base cycle.  In particular,
$A_0N_0.C_0=D$ and the real codimension of $A_0N_0.z$ in $D$ is at
most the real dimension of $C_0$ for every $z\in D$.  Since
$B \supset A_0N_0$, it follows that $\cO\cap C_0\not =\emptyset $ and the
desired dimension bound holds.
\hfill $\square$
\medskip

{\sc Invariant incidence varieties.}
Let $B$ be an Iwasawa--Borel subgroup and $S$ a $B$--Schubert variety
with codim$_Z S = \dim _\mathbb C C_0=q$ and $S\cap C_0\not =\emptyset $.
Using the same type of argument as above, we make the following
observation \cite {H}.

\begin {proposition}
The intersection $S\cap D$ is contained in $\cO$,
and the intersection $\cO\cap C_0$ with the base cycle
is nonempty and transversal at each of its points.  Furthermore,
each of its components is an $A_0N_0$--orbit.
\end {proposition}

Note that the above Schubert varieties are completely determined
topologically
by the Poincar\'e dual $PD(C_0)$ in $H_*(Z;\mathbb Z)$. In particular,
there exist such varieties and given one we refer to the components
of $S\cap \cO$ as {\it Schubert slices}.
Using deeper considerations, we have the following improvement
of the above proposition (see \cite {HuW3}).

\begin {proposition}  \label {Schubert slice}
Every Schubert slice $\Sigma $ intersects every cycle
$C\in\Omega _W(D)$ in exactly one point, and that intersection is
transversal.
\end {proposition}

Turning to incidence varieties, we consider $B$ as above and let $Y$ be any
closed $B$--invariant subvariety of $Z$.  Then the incidence variety
$A_Y:=\{ C\in \Omega :C\cap Y\not =\emptyset\}$
is a closed, $B$--invariant algebraic subvariety of $\Omega $.  It
is a proper subvariety if and only if $C_0\not \in Y$.  We 
only consider this case. Since $B$ has only finitely
many orbits in $Z$, there are only finitely many candidates for
$A_Y$\,.

\medskip
Consider the special case $Y=S\setminus \cO$, where
$S$ is a Schubert variety containing a Schubert slice $\Sigma $
as above.  Here, due to Proposition \ref {Schubert slice}, we
refer to $S$ as a {\it transversal Schubert variety}.

Recall that $Y$ has the structure of a very ample
Cartier divisor. Let $\Gamma (S,{\cal O}(*Y))$ denote the space
of meromorphic functions on $S$ with poles only on $Y$ and let
$\Gamma (\Omega ,{\cal O}(*A_Y)$ be the analogously defined space
of functions on $\Omega $.  The {\it trace transform}
${\cal T}_S:\Gamma (S,{\cal O}(*Y))\to \Gamma(\Omega ,{\cal O}(*A_Y))$
is defined by ${\cal T}_S(f)(C)={\sum }_{p\in C\cap S}f(p).$
(\cite {BK},\cite {BM}, \cite [Appendix]{HuS}).  Of
course this is first defined at the cycles $C$ which intersect
$S$ generically.  This resulting function in
$\Gamma (\Omega,{\cal O}(*A_Y))$ arises via analytic continuation.
With care about cancellations in the defining sum,
one proves 

\begin {theorem} \label{is_hypersurface} {\rm \cite {BK}}
Given $C\in A_Y$ there exists $f\in \Gamma (S,{\cal O}(*Y))$
such that the polar set ${\pi}({\cal T}_S(f))$ contains
$C$.  In particular, $A_Y$ is a $B$--invariant complex hypersurface.
\end {theorem}

Incidence varieties in $\Omega $ which are hypersurfaces are denoted
by $H_Y$ and are called {\em incidence hypersurfaces}.  If
$S$ is a transversal Schubert variety and $Y=S\setminus \cO$, then,
by Proposition \ref{codimension}, the incidence hypersurface $H_Y$ is
contained in the complement of $\Omega _W(D)$ in $\Omega $.

\medskip
One of our original goals was to prove that $\Omega _W(D)$ is
a domain of holomorphy in $\Omega $, i.e., that given a divergent
sequence $\{ C_n\} \subset \Omega _W(D)$ there exists a function
$f\in {\cal O}(\Omega _W(D))$ with
$\lim_{n\to \infty} \vert f(C_n\vert = \infty$.  The
following criterion (see, for example, \cite{GuR}), formulated in the 
setting where $\Omega$ is affine, is useful for this.

\begin {theorem} \label {Stein}
If for every point $C$ in the boundary $\text{\rm bd}(\Omega _W(D))$
in $\Omega $ there exist a complex hypersurface $H$ in $\Omega $
which is contained in $\Omega \setminus \Omega _W(D)$ with
$C\in H$, then $\Omega _W(D)$ is a domain of holomorphy.
\end {theorem}

One consequence of our recent work is 
that for every $C\in \text{bd}(\Omega _W(D))$ there exists
an incidence hypersurface $H_Y$ defined by a transversal
Schubert variety (\cite {HuW3}, also Section \ref{sec7}).
\medskip

\noindent {\bf Domains defined by invariant hypersurfaces.}
\medskip

Let $B \subset G$ be an Iwasawa--Borel subgroup and $H$ a $B$--invariant
complex hypersurface in $\Omega $.  The family $\{ g_0H\}_{g_0\in G_0}$
consists of all such hypersurfaces which are equivalent in the sense
that the variation within the family only depends on the choice of $B$.
Since $H$ is $A_0N_0$--invariant for some Iwasawa decomposition
$G_0=K_0A_0N_0$, this family is the same as $\{ k_0H\}_{k_0\in K_0}$.
The connected component of $\bigcap_{k_0\in K_0} (\Omega \setminus k_0H)$
containing the neutral point $x_0\in \Omega $ is a $G_0$--invariant domain
defined by $H$.  We denote it $\Omega _H$\,.
\smallskip

If $S$ is a transversal Schubert variety, we have the associated
incidence hypersurface $H=H_Y$ where $Y=S\setminus \cO$.  The
connected component $\Omega _S(D)$ containing $x_0$ of the intersection
of all $\Omega _H$ where $H=H_Y$ is an incidence hypersurface associated
to such a Schubert variety is referred to as the
{\it Schubert domain associated to $D$ in $\Omega $}.
\smallskip

Finally, the {\it Iwasawa domain} is defined to be the connected component
containing $x_0$ of the intersection of all the $\Omega _H$, i.e.,
as $H$ ranges over all the (finitely many) $B$--invariant hypersurfaces.

\begin {proposition}
$\Omega _H$, $\Omega _S(D)$ and $\Omega _I$ are $G_0$--invariant
domains of holomorphy in $\Omega $.
\end {proposition}

\noindent {\bf Proof.}
By definition these domains are $G_0$--invariant, connected open subsets
of $\Omega $. Every boundary point of such a domain is contained in a
hypersurface $k_0H$ which is contained in its complement in $\Omega $.
The desired result then follows from Theorem \ref {Stein}
\hfill $\square$

From the definitions, $\Omega _I\subset \Omega _S(D)$ and, if $H$ is a Schubert
incidence hypersurface then
$\Omega _S(D)\subset \Omega _H$\,.

\medskip
The domain $\Omega _I$ can be viewed from several different perspectives.
For example, it is the same as the {\it polar} $\widehat X_0$ which
is defined to be
$\{ gx_0\in \Omega :g\in G \text{ and }X_0\subset g\kappa _0\} $
where $X_0$ is the closed $G_0$--orbit in $\widetilde Z=G/B$ and
$\kappa _0$ is the open $K$--orbit in $\widetilde Z$ (see \cite {Ba}).
R. Zierau remarked that ${\widehat X_0}=\Omega _I$ (see \cite {HuW3}).
The polar can also be regarded as a type of cycle space
(see \cite {GM}).  
\medskip

\noindent {\bf The equality $\Omega _I=\Omega _{AG}$\,.}
\medskip

As noted above, the domain $\Omega _I$ has various guises.  
From the point of view of holomorphic extension of certain
special functions on $\Omega_0$, it is shown in \cite {KS} that
$\Omega _{AG}\subset \Omega _I$ for the classical groups.  From
the polar viewpoint it is shown in \cite{Ba} that
$\Omega_I \subset \Omega_{AG}$\,.   Equality of these domains
comes down to the opposite inclusion, $\Omega _{AG}\subset \Omega _I$\,.
That was proved in \cite{H} using the plurisubharmonic function $\rho$
that is defined by the adapted complex  structure; see Proposition
\ref {easy inclusion}.  Now we put all this together.

\begin {theorem} \label {I=AG}
$\Omega _I=\Omega _{AG}$\,.
\end {theorem}

\noindent {\bf Proof.}
For the inclusion $\Omega _{AG}\subset \Omega _I$ let
$\rho $ be the $G_0$--invariant strictly plurisubharmonic function
of Proposition \ref {easy inclusion}.  Suppose to the contrary that some
hypersurface $H$ which is invariant under an Iwasawa--Borel subgroup
$B$ has nonempty intersection with $\Omega _{AG}$ or, equivalently,
for some Iwasawa decomposition $G_0=K_0A_0N_0$ and some
$x_1\in \Omega _{AG}$ the orbit $AN.x_1$ is not open.
\smallskip

The orbit $A_0N_0.x_0=\Omega _0$ is totally real of half dimension
in $\Omega $.  On the other hand, no $A_0N_0$--orbit in
$AN.x_1\cap \Omega _{AG}$ is of this type.  Now $\rho (x_0)=0$
and $\rho >0$ along all other $G_0$--orbits.  For $x\in \Sigma _0$
near $x_0$, i.e., for $\rho (x)>0$ sufficiently small, $A_0N_0.x$
must still be totally real.
\smallskip

Let $r$ be the smallest value of $\rho \vert \Sigma _0$ such that
for some $x_1\in \{ \rho =r\} $ the orbit $A_0N_0.x_1$ is not
not totally real. Such a value must exist, because the total reality
of an $A_0N_0$--orbit is equivalent to the openness of the corresponding
$AN$--orbit; in particular, we might as well let $x_1$ be the same
point as that which was denoted by $x_1$ at the outset.
\smallskip

Apply Proposition \ref {psh} to the complex manifold
$X=AN.x_1\cap \Omega _{AG}$, the strictly plurisubharmonic function
$\rho|_X$ and the real submanifold $M:=A_0N_0.x_1$. It follows
that $\dim_\mathbb RM\leqq \dim_\mathbb C X$. But this implies that
$\dim_\mathbb R M\leqq \dim_\mathbb R \Omega _0-2$, contrary to the fact
that the $A_0N_0$--orbits in $\Omega _{AG}$ all have the same dimension
as $\Omega_0$\,.
\smallskip

Our proof of the inclusion $\Omega _I\subset \Omega _{AG}$ follows
that in \cite {Ba}.  In the flavor of the present paper we use the
fact that the $G_0$--action on $\Omega _I$ is proper.  That follows
from the Kobayashi hyperbolicity of $\Omega _I$
(\cite {H}, see Section \ref{sec9}).

\smallskip
Suppose that there is a sequence
$\{ z_n\}\subset \Omega _{AG}\cap \Omega _I$ with
$z_n\to z\in \text{bd}(\Omega _{AG})\cap \Omega _I$\,.
From the definition of $\Omega _{AG}$, it follows that
there exist $\{ g_m\}\subset G$ and
$\{ x_m\} \subset \exp(i\omega _0)$ such that $g_m(x_m)=z_m$\,.
Write $g_m=k_ma_mn_m$ in a $K_0A_0N_0$ decomposition of $G_0$\,. Since
$\{ k_m\} $ is contained in the compact group $K_0$\,, it may be
assumed that $k_m\to k$; therefore that $g_m=a_mn_m$\,.
Since $\omega _0$ is relatively compact in ${\mathfrak a}$, it
may also be assumed that $x_m\to x\in \text{\rm c}\ell(\exp(i\omega _0))$.
Thus $x_m=T_mx_0$, where $\{ T_m\} \subset \exp(i\omega _0)$
and $T_m\to T$.  Write $a_mn_n(x_m)=a_mn_mT_mx_0=
\tilde a_m\tilde n_mx_0$, where $\tilde a_m=a_mA_m$
and $\tilde n_m=T_m^{-1}n_mT_m$ are elements of $A$ and
$N$, respectively.
Now $\{ z_m\}$ and the limit $z$ are contained in $\Omega _I$
which is in turn contained in $AN\cdot x_0$.
Furthermore, $AN$ acts freely on this orbit.
Thus $\tilde a_m\to \tilde a\in A$ and
$\tilde n_m\to \tilde n\in N$ with
$\tilde a\tilde n\cdot x_0=z$.  Since $T_m\to T$, it follows that
$a_m\to a\in A_0$ and $n_m\to n\in N_0$ with $an.x=z$.
Since $z\not\in \Omega _{AG}$, it follows that
$x\in \text{\rm bd}(\exp(i\omega _0))$, and $z\in \Omega _I$
implies that $x\in \Omega _I$\,.

\smallskip
On the other hand, since $x\in \text{bd}(\exp(i\omega _0))$, the
isotropy group $G_x$ is noncompact.  But $\Omega _I$ is
Kobayashi hyperbolic (\cite {H}, see Section \ref{sec9}).  Therefore
the $G$--action on $\Omega _I$ is proper (see e.g. \cite {H}) and
consequently $x\notin \Omega _I$, which is a contradiction.
That completes the proof.
\hfill $\square$

\section{Transversal Schubert Varieties.} \label{sec7}
\setcounter{equation}{0}

In this Section we outline the methods introduced in \cite {HuW3}
and their main applications, in particular the fact that
$\Omega _W(D)=\Omega _S(D)$.  The intermediate results, which follow
from a sort of triality, are of complex analytic interest.  For
example, at every boundary point $p\in \text{bd}(D)$ we construct a
$(q+1)$-codimensional Schubert variety $S$ with
$p\in S$ and $S\subset Z\setminus D$.
Due to the existence of the family $\Omega_W(D)$ of $q$-dimensional cycles in $D$, 
this exhibits the maximal possible degree of holomorphic convexity.
It should be underlined that $S$ arises as an
extension of a complex analytic manifold in the the boundary orbit
and therefore its relation to the signature of the Hessian of a boundary 
defining function is unclear.
\medskip

The Borel groups $B$ considered here are Iwasawa--Borel subgroups,
the ones that contain an Iwasawa component $A_0N_0$\,.  The
Schubert varieties are always those which are closures $S$ of orbits
$\cO$ in $Z$ of Iwasawa--Borel subgroups, and $Y:=S \setminus \cO$.  If
codim$_Z S = \dim_\mathbb C C_0=q$ and $S\cap C_0\not =\emptyset $, then
$S$ is called a transversal Schubert variety (see Proposition 
\ref {Schubert slice}).
\smallskip

As usual $D$ is an open $G_0$-orbit in $Z=G/Q$ and
$C_0=K_0.z_0$ is the base cycle.  Here we view $\Omega _W(D) \subset \Omega$ 
with $\Omega := G.C_0$ in the cycle space ${\cal C}^q(Z)$.
In other words, in the case when $G^0_{C_0}=K$ we do not replace
$\Omega $ by the finite cover $G/K$ even if the $G$--stabilizer of $C_0$
properly contains $K$, and we do not exclude the hermitian holomorphic
case where $\Omega =G.C_0$ is an associated compact hermitian symmetric
space.
\smallskip

\noindent {\bf Duality.}
\medskip

Throughout this subsection
$\gamma \in \text{\rm Orb}_ZG_0$ (resp. $\kappa \in \text{\rm Orb}_ZK$) denotes
one of the finitely many $G_0$-orbits (resp. $K$-orbits) in
$Z$.  The first example of duality was proved at the level of
open $G_0$-orbits \cite {W2}: {\it Every open $G_0$-orbit $\gamma $
contains a unique compact $K$-orbit $\kappa $}.

Let us reformulate this in a way that makes sense for all
$G_0$- and $K$-orbits.  For this first observe that, since
$G_0\cap K=K_0$, the intersection $\gamma \cap \kappa $ of any
two $G_0$- and $K$-orbits is $K_0$-invariant.

For $z\in \gamma \cap \kappa $ we refer to the orbit $K_0.z$
as {\it isolated} if it has a neighborhood $U$ in $\gamma $
so that $\kappa \cap U=K_0.z$. Finally, let us say that
$(\gamma ,\kappa )$ is a dual pair if $\gamma \cap \kappa $ contains
an isolated $K_0$-orbit.
\smallskip

Matsuki's duality theorem, which is an extension of the above
statement for open orbits, states that there is a 
bijective map $\mu :\text{\rm Orb}_ZG_0\to \text{\rm Orb}_ZK$ such that
$(\gamma ,\kappa )$ is a dual pair if and only if
$\kappa =\mu (\gamma )$ \cite {M}.
\smallskip

In the original version, duality was defined by $\gamma \cap \kappa $
being nonempty and compact, but it was soon realized that this is
equivalent to the condition that $\gamma \cap \kappa $ be a single
$K_0$--orbit.  The more recent proofs, which
involve the Morse theory related to a certain moment map
(\cite {MUV}, \cite {BL}), use the above weaker notion.  However, in the
end, {\it if an intersection $\gamma \cap \kappa $ contains an isolated
$K_0$--orbit, then it is a $K_0$--orbit}, and 
the intersection $\gamma \cap \kappa $ is transversal along that
orbit.
\smallskip

In our work we use the following {\it non-isolation property}
which is implicit in the proofs in \cite {MUV} and \cite {BL}.

\begin {proposition} \label {non-isolation}
If $(\gamma ,\kappa )$ is not a dual pair, then every $K_0$--orbit 
$K_0.z$ in $\gamma \cap \kappa $ is contained in a $K_0$--invariant
locally closed submanifold $M\subset \gamma \cap \kappa $ with
$\dim_\mathbb RM=\dim_\mathbb RK_0.z+1$.
\end {proposition}

This symplectic approach also yields
information about the topology of the $G_0$-orbits \cite {HuW3}:

\begin {proposition}
If $(\gamma ,\kappa )$ is a dual pair, but $\gamma \cap \kappa \ne \emptyset$,
then the $K_0$-orbit
$\gamma \cap \kappa $ is a $K_0$--equivariant strong deformation retract
of $\gamma $.
\end {proposition}

\noindent {\bf Triality.}
\medskip

Here we describe a sort of triality where, in addition to the orbits
$\gamma $ and $\kappa $ above, we incorporate Schubert
varieties of Iwasawa Borel subgroups $B$.  For this we fix
an Iwasawa decomposition $G_0=K_0A_0N_0$ and $B$ containing
$A_0N_0$.
For $\kappa \in \text{\rm Orb}_ZK$ let $\text{c}\ell (\kappa )$ denote
its closure in $Z$ and define ${\cal S}_\kappa $ to be the
set of all $B$--Schubert varieties $S$ such that
codim$_ZS= \dim_\mathbb C\kappa $ and
$S\cap \text{c}\ell(\kappa )\not =\emptyset $.
\smallskip

The Schubert varieties of a fixed Borel subgroup generate the
integral homology of $Z$ and consequently ${\cal S}_\kappa $
is determined by the topological class of $\text{c}\ell (\kappa )$;
in particular, it is non-empty.
\smallskip

The following can be regarded as a statement of {\it triality}.

\begin {theorem} \label {triality}
If $(\gamma ,\kappa )$ is a dual pair, then the following hold
for every $S\in {\cal S}_\kappa $.
\begin {enumerate}
\item $S \cap \text{\rm \,c}\ell(\kappa)$ is contained in
$\gamma \cap \kappa $ and is finite.
If $x \in S \cap \kappa $, then
$(AN)(x) = B(x) = \cO$, where $S = \text{\rm \,c}\ell(\cO)$, and $S$ is
transversal to $\kappa $ at
$x$ in the sense that the real tangent spaces satisfy
$T_{x}(S) \oplus T_{x}(\kappa) = T_{x}(Z)$.
\item The set $\Sigma = \Sigma(\gamma,S,x) := A_0N_0(x)$ is
open in $S$ and closed in $\gamma$; in particular it is a locally
closed complex submanifold of $Z$.
\item Let {\rm c}$\ell(\Sigma)$ and {\rm c}$\ell(\gamma)$
denote closures in $Z$.  Then the map
$K_0 \times \text{\rm \,c}\ell(\Sigma) \to \text{\rm \,c}\ell(\gamma)$,
given by $(k,z) \mapsto k(z)$, is surjective.
\end {enumerate}
\end {theorem}

\noindent {\bf Proof.}
Let $x \in S \cap \text{\rm \,c}\ell(\kappa)$.
Since ${\mathfrak g} = {\mathfrak k} + {\mathfrak a}  + {\mathfrak n}$
is the complexification of the Lie algebra
version ${\mathfrak g}_0 = {\mathfrak k}_0 + {\mathfrak a}_0 +
{\mathfrak n}_0$ of $G_0 = K_0A_0N_0$\,,
we have $T_x(AN(x)) + T_x(K(x)) = T_x(Z)$.
As $x \in S = \text{\rm \,c}\ell(\cO)$ and $AN \subset B$,
we have $\dim AN(x) \leqq \dim B(x) \leqq \dim \cO = \dim S$.
Furthermore $x \in \text{\rm \,c}\ell(\kappa)$. Thus
$\dim K(x) \leqq \dim \kappa$.
If $x$ were not in $ \kappa$, this inequality would be strict, in violation
of the above additivity of the dimensions of the tangent spaces.
Thus $x\in \kappa $ and $T_x(S) + T_x(\kappa) = T_x(Z)$.
Since $\dim S + \dim \kappa = \dim Z$ this sum is direct, i.e.,
$T_x(S) \oplus T_x(\kappa) = T_x(Z)$.  Now also
$\dim AN(x) = \dim S$ and
$\dim K(x) = \dim \kappa$.  Thus $AN(x)$ is open in $S$, forcing
$AN(x) = B(x) = \cO$.  We have already seen that $K(x)$ is open
in $\kappa $, forcing
$K(x) = \kappa $.  For assertion 1 it remains only to show that
$S \cap \kappa $ is contained in
$\gamma $ and is finite.
\smallskip

Denote $\widehat{\gamma} = G_0(x)$.
If $\widehat{\gamma} \ne \gamma$, then $(\widehat{\gamma},\kappa)$ is
not dual, but $\widehat{\gamma} \cap \kappa$ is nonempty because it contains
$x$.  By the non--isolation property (\ref {non-isolation}), we have a
locally closed $K_0$--invariant manifold
$M \subset \widehat{\gamma} \cap \kappa$ such that
$\dim M = \dim K_0(x) +1$.  We know $T_x(S) \oplus T_x(\kappa) = T_x(Z)$,
and $K(x) = \kappa$, so $T_x(A_0N_0(x))\cap T_x(M) = 0$.  Thus
$T_x(A_0N_0(x)) + T_x(K_0(x))$ has codimension $1$ in the subspace
$T_x(A_0N_0(x)) + T_x(M)$ of $T_x(\widehat{\gamma})$, which contradicts
$G_0 = K_0A_0N_0$\,.  We have proved that $(S \cap \text{\rm \,c}\ell(\kappa))
\subset \gamma$.  Since that intersection is transversal at $x$, it is finite.
This completes the proof of assertion 1.
\medskip

We have seen that $T_x(AN(x)) \oplus T_x(K(x)) = T_x(Z)$, so
$T_x(A_0N_0(x)) \oplus T_x(K_0(x)) = T_x(\gamma)$, and $G_0(x) = \gamma$.
From the basic properties of a dual part, in particular the
transversality of the intersection $\gamma \cap \kappa $, we have
$\dim A_0N_0(x) = \dim T_x(\gamma) - \dim T_x(\kappa \cap \gamma)
= \dim T_x(Z) - \dim T_x(\kappa) = \dim AN(x) = \dim S$.
Now $A_0N_0(x)$ is open in $S$.
\smallskip

Every $A_0N_0$--orbit in $\gamma$ meets $K_0(x)$ because $\gamma = G_0(x)
= A_0N_0K_0(x)$.  By the transversality of the
intersection $\gamma \cap \kappa $, every such $A_0N_0$--orbit
has dimension at least that of $\Sigma = A_0N_0(x)$.   Since the orbits
on the boundary of $\Sigma$ in $\gamma$ would necessarily be smaller,
it follows that $\Sigma$ is closed in $\gamma$.
This completes the proof of assertion 2.
\medskip

The map $K_0 \times \Sigma \to \gamma$, by $(k,z) \mapsto k(z)$, is surjective
because $K_0A_0N_0(x) = \gamma$.  Since $K_0$ is compact and $\gamma$
is dense in $\text{\rm \,c}\ell(\gamma)$, assertion 3 follows.
\hfill $\square$
\medskip

Let us now indicate the construction given in \cite {HuW3} of
the {\it supporting Schubert variety} at each boundary point
$p\in \text{bd}(D)$.  For this it is convenient to refer
to a $G_0$-orbit in $\text{bd}(D)$ as {\it generic} if
it is open in $\text{bd}(D)$. The transversality property \cite {HuW3}
also gives us

\begin {lemma}
If $\gamma $ is a generic orbit in {\rm bd}$(D)$, $\kappa $
is dual to $\gamma $ and $S=\text{c}\ell (\cO)\in {\cal S}_\kappa $,
then {\rm codim}$_ZS> \dim_\mathbb CC_0$.
\end {lemma}

There are increasing sequences $\{ \cO_k\} $ of $B$-orbits
with $\cO_0=\cO$, $\dim\,\cO_k=\dim\,\cO+k$ and
$\cO_k\subset \text{c}\ell(\cO_{k+1})$, for $k \leqq n - \dim_\C \cO$,
so the above Schubert variety
can be enlarged to obtain the following consequence of
triality.

\begin {theorem}
For every $p\in \text{\rm bd}(D)$ there exists an Iwasawa-Borel
subgroup $B \subset G$ and a $B$--Schubert variety $S$ with
{\rm (1)} {\rm codim}$_ZS=\dim_\mathbb CC_0+1$ and
{\rm (2)} $p\in S$ and $S\subset Z\setminus D$.
\end {theorem}

\noindent {\bf Proof.}
Given a point $p$ in a generic orbit $\gamma $, we may take an
appropriate conjugate of the Iwasawa-Borel subgroup of
Theorem \ref {triality} to obtain $S$ satisfying all of the required conditions
except that it may be too small.  In that case we enlarge
it to be $(q+1)$-codimensional by the above procedure. By
Proposition \ref {codimension} this Schubert variety is also contained
in $Z\setminus D$.
\smallskip

If $p\in \text{bd}(D)$ is not generic, then nevertheless it is
the limit $p_n\to p$ with $\gamma =G_0.p_n$ generic and with
$B_n$-Schubert varieties $S_n$ which have the desired properties.
\smallskip

Since $\{ S_n\}$ is contained in a compact space of cycles and
$\{ B_n\} $ can be regarded as a sequence in the (compact)
minimal $G_0$-orbit, by going to subsequences, $S_n\to S$ and
$B_n\to B$, where $S$ is a $B$-Schubert variety with the
desired properties. \hfill $\square$
\medskip

\noindent {\bf The equality $\Omega _W(D)=\Omega _S(D)$\,.}
\medskip

This equality is equivalent to the following statement.

\begin {theorem} \label{trans}
Given $C\in \text{\rm bd}(\Omega _W(D))$ there exists a transversal
Schubert variety $S=\cO \cup \,Y$ with $C\cap Y\not=\emptyset$.
\end {theorem}

We outline a constructive proof of the existence of
such an $S$. See \cite {HuW3} for the details.
\smallskip

For any initial $p\in C\cap \text{bd}(D)$ we begin with a smallest
$B$--Schubert variety $S$ such that (i)  $S\cap D\not=\emptyset$, (ii)
$p$ is in a component of $Y=S\setminus \cO$ that is contained
in the complement $Z\setminus D$, and (iii) all components of
$Y$ which contain $p$ are themselves in $Z\setminus D$.
Define $\delta \geqq 0$ by
$\text{ \rm codim}_ZS=\dim_\mathbb C+1-\delta$.
If $\delta =0$, then there is nothing to prove.
\smallskip

Now let $\delta >0$.  Let $U$ denote the union of the components of $Y$ which
have empty intersection with $D$ and $V$ be the union of the remaining
components.  $V$ is non-empty by the minimality assumption on $S$.

Using intersection theory information, roughly speaking the fact that the 
intersection $C_0.U$ is zero in homology, one proves that $S$ can be replaced by a
lower--dimensional $B$--Schubert variety $S_1$ which is contained
in $V$ with $p$ being replaced by $p_1\in C\cap \text{bd}(D)$.

Continuing in this way, one eventually determines a transversal
Schubert variety $S_\delta $ and a point
$p_\delta \in C\cap \text{bd}(D)$ which is also contained
in $Y_\delta $.\hfill $\square $

\begin {corollary} \label {W=S}
$\Omega _W(D)=\Omega _S(D)$
\end {corollary}

\noindent {\bf Proof.}
Given $C\in \text{bd}(\Omega _W(D))$, the above result yields an
incidence hypersurface $H_Y$ which contains $C$
and which is contained in the complement
$\Omega \setminus \Omega _W(D)$; see Theorem \ref{is_hypersurface}.
\hfill $\square$

\begin {corollary}
The cycle space $\Omega _W(D)$ is a Stein domain in $\Omega $.
\end {corollary}

\noindent {\bf Proof.}
When $\Omega $ is affine, this follows immediately
from Theorem \ref {Stein} and the equality $\Omega _W(D)=\Omega _S(D)$.
In the hermitian holomorphic cases, where $\Omega$ is the associated
compact symmetric space, note that, since Pic$(\Omega )\cong \mathbb Z$,
it follows that the complement of any $H_Y$ in $\Omega $ is affine.
Since $\Omega_W(D)\subset (\Omega \setminus H_Y)$, the result follows 
for the same reasons as above.
\hfill $\square$
\medskip

\noindent {\bf Spaces of cycles in lower dimensional $G_0$--orbits.}
\medskip

We recall the setting of \cite {GM}.  For $Z=G/Q$,
$\gamma \in \text{\rm Orb}_Z(G_0)$ and $\kappa \in \text{\rm Orb}_Z(K)$ its
dual, let $G\{\gamma\}$ be the connected component of the
identity of
$\{ g\in G:g(\kappa )\cap \gamma \ \text{is non-empty and compact}\ \}$.
Note that $G\{\gamma\}$ is
an open $K$-invariant subset of $G$ that contains the identity.
Define ${\cal C}\{\gamma\}:=G\{\gamma\}/K$.  Finally, define
${\cal C}$ as the intersection of all such cycle spaces ${\cal C}\{\gamma\}$
as $\gamma $ ranges over $\text{\rm Orb}_Z(G_0)$ and $Q$ ranges over all
parabolic subgroups of $G$.

\begin {theorem}
${\cal C}=\Omega _{AG}$\,.
\end {theorem}

This result was checked in \cite {GM} for classical and hermitian exceptional
groups using case by case computations, and the authors of \cite{GM}
conjectured it in general.  As will be shown here, it is a
consequence of the fact that
$\Omega _W(D)=\Omega _S(D)$ in the special case where
$D$ is an open $G_0$-orbit in $G/B$
and of the following general result \cite[Proposition 8.1]{GM}.

\begin {proposition}
$\left ( \bigcap_{D \subset G/B \text{ open}}\ \Omega_W(D) \right )
\subset {\cal C}$.
\end {proposition}

\noindent {\bf Proof of Theorem.}
The polar $\widehat {X_0}$ in $Z=G/B$ coincides with the
cycle space ${\cal C}_Z(\gamma _0)$, where $\gamma _0$ is the unique
closed $G_0$-orbit in $Z$. As was shown above, this agrees
with the Iwasawa domain $\Omega _I$\, which in turn is contained
in every Schubert domain $\Omega _S(D)$. Thus, for every open
$G_0$-orbit $D_0$ in $Z=G/B$ we have the inclusions

\begin {equation*}
\left ( {\bigcap}_{D \subset G/B \text{ open}}\ \Omega_W(D) \right )
\subset {\cal C}\subset {\cal C}_Z(\gamma _0)=
\widehat{X_0} =\Omega _I\subset \Omega _S(D_0)=\Omega _W(D_0).
\end {equation*}
Intersecting over all open $G_0$-orbits $D$ in $G/B$, the equalities
\begin {equation*}
\left ( {\bigcap}_{D \subset G/B \text{ open}}\ \Omega_W(D) \right )
={\cal C}=\Omega _I=
\left ( {\bigcap}_{D \subset G/B \text{ open}}\ \Omega_W(D) \right )
\end {equation*}
are forced, and ${\cal C}=\Omega _{AG}$ is a consequence of
$\Omega _I=\Omega _{AG}$.\hfill $\square $

\section{Cycle Domains in the Hermitian Case.} \label{sec8}
\setcounter{equation}{0}

In this section $G_0$ is a group of hermitian type, in other words
$\cB = G_0/K_0$ is a bounded symmetric domain.  We give a concrete 
description of $\Omega _W(D)$.  As indicated
in Section \ref {sec2}, either $\Omega _W(D)=\cB$ or $\overline{\cB}$, or 
$G.C_0=\Omega $
is affine.  Thus it is enough to consider the latter case.
At first we will replace $\Omega _W(D)$ by the connected component
containing $x_0$ in its preimage in $G/K$ and denote the latter
by $\Omega $.
\medskip

Choosing a system of roots in the usual way, we regard the 
bounded symmetric domain ${\cal B}$ of $G_0$ as the $G_0$-orbit of the
neutral point $x_0\in X=G/P_-$ and its complex conjugate as the
orbit of the analogous point $\overline{x_0} \in \overline{X} =G/P_+$.  We view
$z_0=(x_0,\overline{x_0})$ as the base point for
$\Omega =G.z_0=G/K\hookrightarrow X\times \overline{X}$, i.e., the open
$G$-orbit by its diagonal action.  We identify
${\cal B}\times \overline {\cal B}$ with its image under the natural
embedding ${\cal B}\times \overline {\cB}\hookrightarrow X\times \overline{X}$
and note that this lies in $\Omega $.
\smallskip

The following is proved by a reduction to the polydisc case \cite {BHH}.

\begin {proposition}
$\Omega _{AG}={\cal B}\times \overline {\cal B}$
\end {proposition}

In \cite {WZ1} it was shown that
$\Omega _W(D)\subset {\cal B}\times \overline {\cal B}$.  Thus
Corollary \ref {W=S} together with Theorem \ref {I=AG} imply the
following characterization.
\begin {theorem} \label{herm_case} {\rm (\cite {HuW3}, \cite {WZ3})}
If $G_0$ is of hermitian type, then either {\rm (i)} $\Omega$ is
the compact dual symmetric space to $G_0/K_0$ and $\Omega _W(D)$ is ${\cal B}$
or $\overline {\cal B}$, or {\rm (ii)} $\Omega =G/\tilde K$ is affine with 
$\tilde K/K$ finite and $\Omega _W(D)={\cal B}\times \bar {\cal B}$.
In case {\rm (ii)} $\Omega _W(D)$ lifts bijectively to $G/K$.
\end {theorem}

\noindent {\bf Proof.}
It is enough to consider the case where $\Omega $ is affine.
As we have seen above, after lifting to $\Omega =G/K$,
$\Omega _W(D)\subset {\cal B}\times \overline {\cal B}=\Omega _{AG}$.
But, by Theorem \ref {I=AG}, $\Omega _{AG}=\Omega _I$\,, and
$\Omega _W(D)=\Omega _S(D)$ by Corollary \ref {W=S}. Since by definition
$\Omega _I\subset \Omega _S(D)$, the result follows at the level of
$G/K$.  Furthermore, ${\cal B}\times \overline {\cal B}$ is a cell.
Thus it agrees with its image in $G/\tilde K$ by the canonical
finite covering map (see Proposition \ref {cover_trivial}).
\hfill $\square$

\begin {remark}{\em
This theorem can also be proved using results from \cite {GM}.}
\end {remark}

\section{Kobayashi Hyperbolicity.} \label{sec9}
\setcounter{equation}{0}

The Kobayashi pseudometric $d_K$ is defined on any complex
manifold $X$; see Section \ref {sec5}.  If it is a metric, i.e.,
$d_K(x,y)>0$ if $x\not =y$, then $X$ is said to be Kobayashi
hyperbolic.  One checks that it is a metric on bounded domains
and vanishes identically in the case where $X$, e.g., is the complex plane.
Since holomorphic maps are Kobayashi-distance decreasing, if there
exists a non-constant holomorphic map $f:\mathbb C\to X$, then
$X$ is not Kobayashi-hyperbolic.
\smallskip

With this in mind it is useful to introduce a weaker notion of
hyperbolicity: $X$ is said to be {\it Brody hyperbolic} if there
are no non-constant holomorphic maps $f:\mathbb C\to X$.
\smallskip

If $X$ is Kobayashi hyperbolic, then its group Aut$_{\cal O}(X)$
of holomorphic automorphisms, equipped with the compact--open
topology, is a Lie group acting properly on $X$.
This group is in general not semisimple, but nevertheless
semisimple subgroups act properly.

\begin {proposition}
Let $G$ be a semisimple Lie group with finitely many components such that
$G^0$ has finite center.  Then any smooth almost effective action of 
$G$ by holomorphic transformations on a hyperbolic manifold $X$ is proper.
\end {proposition}

\noindent {\bf Proof.}
Let $L:= \text{Aut}(X)^0$.  By hypothesis on $G$ and its action
on $X$, the action has finite kernel, so for the proof we may assume
that $G$ is connected and contained in $L$.

Now let $R$ be the radical of $L$ and $L=R\cdot S$ be a Levi--Malcev
decomposition, where $S$ is semisimple.  The corresponding Lie algebra
decomposition is $\gl =\mathfrak r\rtimes \mathfrak s$.
Let $\Gamma :=R\cap S$.  It is a discrete normal (thus central) subgroup of $S$.
Projection $L\to L/R=S/\Gamma =:\bar S$, followed by the adjoint 
representation of $\bar S$, is a  Lie group morphism
$\varphi :L\to {\text Ad}(\bar S) \subset GL(\gs)$, and 
Ker$(\varphi|_G)$ is a discrete
central subgroup of $G$, thus finite.  By construction $\varphi(G)$ is
a semisimple subgroup of $GL(\gs)$, thus closed in $GL(\gs)$.
As Ker$(\varphi|_G)$ is finite now $G$ is closed in $L$.  Since the
action of $L$ on $X$ is proper, it follows that the action of $G$ is also 
proper.
\hfill $\square$
\medskip

\noindent {\bf Families of hypersurfaces.}
\medskip

It is a classical result that the complement of the union
of $(2m+1)$-hyperplanes in general position in $\mathbb C\mathbb P_m $
is Kobayashi hyperbolic \cite {D}.  Let us make the notion
of {\it general position} precise in a context which is appropriate
for our applications.
\smallskip

Since the complex manifolds which we consider are embedded in
projective spaces by sections of line bundles, it is natural
to regard a ``point'' as being in the projectivization
$\mathbb P(V^*)$ of the dual space of a complex vector space
and a ``hyperplane'' as a point in $\mathbb P(V)$.
We regard a subset $S\subset \mathbb P(V)$ as parameterizing
a family of hyperplanes in $\mathbb P(V^*)$.
A non--empty subset $S\subset \mathbb P(V)$ is said to have the {\it normal
crossing property} if for every $k\in \mathbb N$ there exist
$H_1,\ldots H_k\in S$ so that for every subset
$I\subset \{ 1,\ldots ,k\} $  the intersection
$\bigcap _{i\in I}H_i$ is $|I| $-codimensional.  If
$|I| \geqq \dim_{\mathbb C}V$, this means that the intersection is empty.
\smallskip

In the sequel $\langle S\rangle$ denotes the complex linear span of $S$ in
$\mathbb P(V)$, i.e., the smallest complex subspace in $\mathbb P(V)$ containing $S$.
If  $\langle S\rangle =\mathbb P(V)$ we say that $S$ is a generating set.

\begin {proposition}
A locally closed, irreducible real analytic subset $S$ with
$\langle S\rangle=\mathbb P(V)$ has the normal crossing property.
\end {proposition}

\noindent {\bf Proof.}
We proceed by induction over $k$. For $k=1$ there is nothing to prove.
Given a set $\{ H_{s_1},\ldots ,H_{s_k}\} $ of hyperplanes
with the normal crossing property and a subset
$I\subset \{s_1,\ldots ,s_k\} $, define
$$
\Delta _I:={\bigcap}_{s\in I}H_{s},
\qquad {\cal H}(I):=\{ s\in S : H_s\supset \Delta _I \}
\qquad {\cal C}\ell_k:=
{\bigcup}_{J\subset \{s_1,..,s_k\},
\Delta _J\ne\emptyset} \ {\cal H}(J).
$$
We wish to prove that $S\smallsetminus {\cal C}\ell _k \ne\emptyset $.
For this, note that each ${\cal H}(I)$ is a real analytic subvariety
of $S$. Hence,
if $S={\cal C}\ell_k$, then $S={\cal H}(J)$ for some $J$ with
$\Delta _J \ne \emptyset $. However,
$\{H\in \mathbb P(V^*):H\supset \Delta _J\} $ is a proper,
linear plane ${\cal L}(J)$
of $\mathbb P(V)$.  Consequently, $S\subset {\cal L}(J)$, and this would
contradict $\langle S\rangle=\mathbb P(V)$.  Therefore, there exists
$s \in S\smallsetminus {\cal C}\ell_k$, or equivalently,
$\{ H_{s_1},\ldots ,H_{s_k},H_s\} $ has the normal crossing property.
\hfill $\square$
\medskip

In the case of finitely many hyperplanes in $\mathbb C\mathbb P_m$
the condition that $H_1,\ldots, H_{2m+1}$ are in general position
is equivalent to their having the normal crossing property, in which case
$\mathbb C\mathbb P(V^*)\smallsetminus \bigcup H_j$ is Kobayashi hyperbolic
(\cite {D}, or see \cite[p. 137]{K}).

\begin {corollary} \label {kobayashi}
If $S$ is a locally closed, irreducible and generating real analytic subset of
$\mathbb P(V)$, then there exist hyperplanes $H_1,\ldots H_{2m+1}\in S$
such that the complement $\mathbb P(V^*)\smallsetminus \bigcup H_j$ is
Kobayashi hyperbolic.
\end {corollary}

Our main application of this result arises in the case where
$S$ is an orbit of the real form at hand.

\begin {corollary}
Let $G$ be a reductive complex Lie group, $G_0$ a real
form, $V^*$ an irreducible $G$-representation space and
$S$ a $G_0$-orbit in $\mathbb P(V)$.  Then there exist
hyperplanes $H_1,\ldots ,H_{2m+1}\in S$ so that
$\mathbb P(V^*)\smallsetminus  \bigcup H_j$ is Kobayashi hyperbolic.
\end {corollary}

\noindent {\bf Proof.}
From the irreducibility of the representation $V^*$, it follows
that $V$ is likewise irreducible and this, along with the identity
principle, implies that $\langle S\rangle=\mathbb P(V)$.
\hfill $\square$
\medskip

\noindent {\bf Invariant hyperbolic domains in $\Omega $.}
\medskip

We begin by briefly discussing the formalism for equivariantly
embedding $\Omega =G/K$ in a complex projective space.  In order
to avoid discussions of infinite-dimensional spaces, we let
$X$ be a $G$-equivariant, smooth, projective algebraic
compactification of $\Omega$.  Since $X$ is rational and we may assume that
$G$ is simply connected, every line bundle $\L\to X$ is a
$G$--bundle; in particular there is a canonically induced action
on the space $\Gamma (X,\L)$ of sections.
\smallskip

Now let $B$ be an Iwasawa-Borel subgroup of $G$ and $H$ be
a $B$-invariant complex hypersurface in $\Omega $.  Since
$H$ is algebraic, it is Zariski open in its closure $\text{c}\ell (H)$
in $X$.
\smallskip

Let $\L\to X$ be the line bundle defined by $\text{c}\ell (H)$ and
$s\in \Gamma (X,L)$ the defining section.  Since $\text{c}\ell (H)$
is $B$--invariant, $s$ is a $B$--eigenvector.
\smallskip

Define $V$ to be the irreducible $G$--representation subspace of
$\Gamma (X,\L)$ that contains $s$.  Let
$\varphi :X\rightharpoonup \mathbb P(V^*)$ denote the canonically
associated $G$--equivariant meromorphic map.

\begin {proposition}
Assume that $G_0$ is not of hermitian holomorphic type.  Then the restriction
$\varphi|_\Omega :\Omega \to \mathbb P(V^*)$ is a $G$--equivariant,
finite--fibered regular morphism with image a quasi-projective
$G$--orbit in $\mathbb P(V^*)$.
\end {proposition}

\noindent {\bf Proof.}
Since $\varphi $ is $G$-equivariant and its set of indeterminacies
is therefore $G$--invariant, the restriction to the open
$G$--orbit is base point free.  The fact that $\varphi|_\Omega $
is finite--fibered is a consequence of $\varphi $ being non-constant,
e.g., $s$ is not $G$--fixed, and the fact that the $G$--isotropy
group $G_{x_0}=K$ is dimension--theoretically a maximal subgroup
of $G$. (It is here that we use the assumption that $G_0$ is not of
hermitian holomorphic type.)
\hfill $\square$

We are now in a position to prove the main theorem of this
section (\cite {H}, \cite {FH}).

\begin {theorem} \label{k_hyper}
If $G_0$ is not of hermitian holomorphic type, $B \subset G$ is an 
Iwasawa--Borel subgroup, and $H$ is a $B$-invariant hypersurface in
$\Omega =G/K$, then $\Omega _H$ is Kobayashi hyperbolic.
\end {theorem}

\noindent {\bf Proof.}
We replace $\varphi $ by its restriction to $\Omega $ and
only discuss that map. By definition every section $\tau
\in V$ is the pull-back
$\varphi ^*(\widetilde \tau )$ of a hyperplane section. Thus,
there is a uniquely defined $B$--hypersurface $\widetilde H$
in $\mathbb P(V^*) $ with $\varphi ^{-1}(\widetilde H)=H$.  Let
$\widetilde \Omega _{\widetilde H}\subset \mathbb P(V^*)$ be defined
analogously to $\Omega _H$, i.e.,  $\widetilde \Omega _{\widetilde H}=
\mathbb P(V^*)\smallsetminus \bigcup _{g\in G_0} g(\widetilde H).$
Applying Corollary \ref{kobayashi} to $\mathbb P(V^*)$ and
$S:=G_0\cdot \widetilde H\subset \mathbb P(V)$,
it follows that the domain
$\widetilde \Omega _{\widetilde H}$ is Kobayashi hyperbolic.
Furthermore, the connected component of
$\varphi ^{-1}(\widetilde \Omega _{\widetilde H})$ which
contains the base point $x_0$ is just the original domain $\Omega _H$.
Since holomorphic maps are distance decreasing and 
$\phi : \Omega_H \to \widetilde \Omega _{\widetilde H}$ is locally
holomorphic, it follows that
$\Omega _H$ is also Kobayashi hyperbolic.
\hfill $\square$
\medskip

Let us summarize what we have presented up to this point.

\begin{summary} \label{summary1} {\em
If $G_0$ is of hermitian holomorphic type, then either
$\Omega$ is the compact hermitian symmetric space dual to the bounded
symmetric domain $\cB = G_0/K_0$\,, and either 
$\Omega _W(D)=\Omega _S(D)={\cal B}$ or
$\Omega _W(D)=\Omega _S(D)=\overline{\cal B}$.
If $G_0$ is of hermitian nonholomorphic type, then $\Omega $ is affine and
\begin {equation*}
\Omega _{AG}=\Omega _I= \Omega_S(D)=\Omega _W(D)={\cal B}
\times \bar {\cal B}.
\end {equation*}
If $G_0$ is not of hermitian type, then, with the usual
convention that $\Omega _W(D)\subset \Omega =G/K$,
\begin {equation*}
\Omega _{AG}=\Omega _I\subset \Omega _S(D)=\Omega _W(D)\subset \Omega _H
\end {equation*}
for $H$ any complex hypersurface which is invariant by an
Iwasawa-Borel subgroup.  All of the domains are $G_0$-invariant,
Stein and Kobayashi hyperbolic and therefore the $G_0$-actions
are proper in every case.}
\end{summary}

\section{The Maximal Domain of Hyperbolicity.} \label{sec10}
\setcounter{equation}{0}

Recall the basic sequence of inclusions
\begin {equation} \label{incl}
\Omega _{AG}=\Omega _I\subset \Omega _W(D)=\Omega _S(D)\subset \Omega _H
\end {equation}
for any complex hypersurface $H\subset \Omega =G/K$ which is
invariant under an Iwasawa-Borel subgroup of $G$, when $D$ is not of
hermitian holomorphic type.  All of these
domains are $G_0$-invariant, Stein and Kobayashi hyperbolic; see
Summary \ref{summary1}.
\medskip

Our goal here is to outline the proof of the following main theorem
of \cite{FH}.

\begin {theorem} \label {main}
The only $G_0$-invariant, Stein, Kobayashi hyperbolic domain in $\Omega $
which contains $\Omega _{AG}$ is $\Omega _{AG}$ itself.
\end {theorem}

This implies that the sequence (\ref{incl}) of inclusions is a sequence of
equalities in the non--hermitian case.  With the
classification in the hermitian case (Section \ref {sec8}), this yields
the following.

\begin {theorem}
In the hermitian holomorphic case, $\Omega _W(D)$ is the bounded domain 
${\cal B}$ or $\overline{\cB}$ associated to $G_0$\,.
In all other cases, $\Omega _{AG}=\Omega _I=\Omega _W(D)=\Omega _S(D)$.
\end {theorem}

The proof of Theorem \ref{main} involves three main steps:
(1) Understanding the invariant theory of the $G_0$--action
on $\Omega $; in particular, the orbit structure on bd$(\Omega _{AG})$.
(2) For every generic boundary point $z$, determining an
${\mathfrak sl}_2$-triple defined over $\mathbb R$ so that the
orbit $S.z=Q_2$ intersects $\Omega _{AG}$ in a $2$--dimensional
polydisk which is the $\Omega _{AG}$ for the group $S$.
(3) Proving Theorem \ref{main} in the case of $S=SL_2(\mathbb C)$.
\medskip

All details which are omitted in the following sketch can be found
in \cite{FH}.  To be consistent with the notation of that
paper, we let $x\in \Omega $ be a general point and $x_0\in \Omega $
the chosen base point with $G_{x_0}=K$.
\medskip

\noindent {\bf The linear model.} 
\medskip

Using a classical linearizing map (see \cite{M1})
we $G_0$--equivariantly embed $\Omega $ in a linear space where
Jordan decomposition can be used in an optimal way.
\smallskip

{\sc The basic map.}
Let $\sigma$ denote complex conjugation of $\gg$ over
${\mathfrak g}_0$\,.  Write $\tau$ for the complex--linear
extension to ${\mathfrak g}$ of the Cartan involution $\theta $ of 
${\mathfrak g}_0$\,, so $\tau$ is the holomorphic involution of 
${\mathfrak g}$ with fixed point set $\gk$.  We extend $\theta$ 
conjugate--linearly to $\gg$, so $\theta$ becomes the conjugate--linear
involution of $\gg$ whose fixed point set is the compact real form
$\gg_u$\,, in other words $\theta$ becomes a Cartan involution of $\gg$.
\smallskip

Define $\eta :G\to \text{Aut}_{\mathbb R}({\mathfrak g})$ by
$\eta (x)=\sigma\circ Ad(x)\circ \tau \circ \text{Ad}(x^{-1})$.  We regard
$\sigma $, $\tau $ and $G$ as operating on Aut$_\mathbb R({\mathfrak g})$
by conjugation, the action of $G$ being via Ad.  Let $N$ denote the normalizer
of $K$ in $G$.  The basic properties of $\eta$ are as follows.
\smallskip

\begin {itemize}
\item  $\eta$ is right $N$--invariant and
defines an embedding $G/N\hookrightarrow
\text{Aut}_\mathbb R({\mathfrak g})$.

\item $\eta$ is equivariant with respect to the left action of $G_0$
on $G$.

\item The image Im$(\eta )$ is closed, $\sigma $--invariant and
$\tau $--invariant, $\sigma (\eta (x))=\eta (x)^{-1}$ and
$\tau (\eta (x))=\eta (\tau (x))$.
\end {itemize}

To see that Im$(\eta )$ is closed we use a
basic lemma of invariant theory (see \cite[p. 117]{Hu}):

\begin {lemma} \label {Humphreys}
Let $V$ be a finite dimensional real vector space,
$H$ a closed reductive algebraic subgroup of $GL _{\mathbb R}(V)$ and
$s\in GL _{\mathbb R}(V)$ an element which normalizes $H$. Regard
$H$ as acting on $GL _{\mathbb R}(V)$ by conjugation. Then,
for a semisimple  $s$  the orbit $H.s $ is closed.
\end {lemma}

To prove that Im$(\eta )$ is closed it is enough to show that
$G.\tau =\{ \text{Ad}(g)\tau \text{Ad}(g)^{-1} : g\in G\} $ is closed in
Aut$_\mathbb C({\mathfrak g})$.  Since $\tau $ is semi-simple
and normalizes $G$ in this representation, this is then immediate.
\medskip

From now on, unless otherwise stated, we view $\eta$ as a
map $\eta :\Omega \to \text{Aut}_\mathbb R({\mathfrak g})$ which is
is essentially a diffeomorphism onto its closed image.  In this
language Lemma \ref {Humphreys} also leads to

\begin {proposition} \label {closed}
If $\eta (x)=s$ is semi-simple, then $G_0.x$ is closed.
\end {proposition}

\noindent {\bf Proof.}
It is enough to show that $G_0.s$ is closed.  By Lemma \ref {Humphreys}
the complex orbit $G.s$ is closed. Define
$\hat \sigma :\text{Aut}_\mathbb R({\mathfrak g})\to 
\text{Aut}_{\mathbb R}({\mathfrak g})$
by $\hat \sigma (\varphi )=(\sigma(\varphi ))^{-1}$.  Here $\sigma $
acts by conjugation as usual.   Note that 
Im$(\eta )$ belongs to the fixed point set Fix$(\hat \sigma )$.  Since
$G.s\cap \text{Fix}(\hat \sigma )$ consists of only finitely many
$G_0$-orbits \cite{Br}, it follows that $G_0.s$ is closed
\hfill $\square$
\smallskip

{\sc Jordan decomposition.}
Let $\eta (x)=u\cdot s$ denote the Jordan decomposition of an
element in Im$(\eta )$.  Then the unipotent factor 
$u\in \text{Aut}_\mathbb R({\mathfrak g})^\circ $,
and $u = \exp(\text{ad}(\nu)) = \text{Ad}(\exp(\nu))$ for some
nilpotent $\nu \in \gg_0$\,.  If $\varphi \in \text{Aut}_\mathbb R({\mathfrak g})$
let $\gg^\varphi$ denote its fixed point set.  The fact $su = us$ can be
expressed $\nu \in \gg^s$.  Also, $\sigma (\eta (x))=\eta (x)^{-1}$ implies
$\eta \in i\gg_0$\,.  Now compute $\eta (\exp(\frac{1}{2}\nu).x)=s$.
In summary we have the following result, describing a lifting
of the Jordan decomposition.

\begin {proposition}   \label {lifting}
For $x\in \Omega $ with Jordan decomposition $\eta (x)=u\cdot s$ there
exists a nilpotent element $\nu \in {\mathfrak g}^s\cap i{\mathfrak g_0}$
such that $u=\text{\rm Ad}(\exp(\nu))$ and 
$\eta (\exp(\frac{1}{2}\nu)\cdot x)=s$.
\end {proposition}

\noindent {\bf Orbit structure.} 
\medskip

Here we outline some basic information on the orbit structure
of the $G_0$-action on $\Omega $.  The main objective is an
understanding of the $G_0$-action on $bd(\Omega _{AG})$.
\medskip

{\sc Closed orbits.}
In any invariant theoretic situation it is of central importance
to move in a systematic way from a point in a non-closed orbit
to a closed orbit in its closure.  Here we accomplish this
by means of special ${\mathfrak sl}_2$-triples.

\medskip
Since $\sigma (s)=s^{-1}$ and $s$ is semisimple, ${\mathfrak g}^s$ is a
$\sigma $-invariant reductive subalgebra of ${\mathfrak g}$.  Let
${\mathfrak g}={\mathfrak h}\oplus {\mathfrak q}$ be its
$\sigma $-decomposition.   Now apply the Jacobson--Morozov Theorem.

\begin {lemma} \label {sltriple}
Let $e\in {\mathfrak g}^s\cap i{\mathfrak g}$
be nonzero and nilpotent.
There exists an ${\mathfrak sl}_2$-triple $(e,h,f)$ in
${\mathfrak g}^{s}$ \text{\rm (}i.e., $[e,f]=h,\ [h,e]=2e$, $[h,f]=-2f$\text{\rm )}
such that  $e,f\in {\mathfrak q}$ and $h\in {\mathfrak h}$.
\end {lemma}

Using this ${\mathfrak sl}_2$-triple we are able to explicitly move
to a closed orbit.

\begin {proposition}
If $\eta(x)=us$ is the Jordan decomposition,
then the orbit $G_0.\eta (x)=G_0.(su)$
contains the closed orbit $G_0.s\subset \text{\rm Im}(\eta )$
in its closure c$\ell(G_0.\eta (x))$.
In particular, $G_0.\eta (x)$
is closed if and only if
$\eta (x)$ is semi-simple.
\end {proposition}

\noindent {\bf Proof.}
Let $u=\text{Ad}(\exp \nu)$ with $\nu$ as in Proposition \ref{lifting}.
Hence, by Lemma \ref{sltriple} there is a ${\mathfrak sl}_2$-triple
$(\nu,h,f)$ ($e=\nu$) such that $[th,\nu]=2t\nu$, i.e., 
Ad$(\exp(th))(\nu)=e^{2t}\nu$
for every $t\in \mathbb R$. Note also that
$\exp (\mathbb R H)\subset G_0\cap \exp ({\mathfrak g}^s)$
by construction of the ${\mathfrak sl}_2$-triple. It follows
that
\begin{eqnarray*}
\eta(\exp (th)\cdot x) & = & \exp (th).(us)=
\text{\rm Ad}(\exp (th))\cdot us\cdot \text{\rm Ad}(\exp(-th))= \\
& = & \text{\rm Ad}(\exp(tH))\cdot u\cdot \text{\rm Ad}(\exp(-tH)) \cdot s= \\
& = & \text{\rm Ad}(\exp(tH))\text{\rm Ad}(\exp \nu) \text{\rm Ad}(\exp (-tH)) \cdot s= \\
& = & \text{\rm Ad}(\exp( e^{2t}\nu)) \cdot s.
\end{eqnarray*}
For $t\to -\infty$ it follows $\exp(tH).(us)=\text{Ad}(\exp(e^{2t}\nu))\cdot s\to s.$
Hence, the closed  orbit $G.s$ lies in the closure of $G.(us).$
In particular $G_0.(us)$ is not closed if $u\ne 1$, i.e., if $\eta(x)$ is
not semisimple. This, together with Proposition \ref{closed} implies that $G.\eta(x)$
is closed if and only if $\eta(x)$ is semi-simple.
Recall that the image Im$(\eta)$ is closed. This forces
$s\in \text{Im}(\eta)$ and the proof is now complete.
\hfill $\square$
\smallskip

{\sc Elliptic elements.} 
An element $x\in \Omega $ is called {\it elliptic} if
$\eta (x)$ is elliptic in the sense that it is semisimple
and all its eigenvalues have absolute value $1$.
Let $\Omega _{ell}$ denote the set of elliptic elements.
\smallskip

If $x \in G_u$\,, so $\theta(x) = x$, then $\theta \eta (x)=\eta (x)\theta$,
so $\eta(x)$ is elliptic.  Therefore $G_u.x_0\subset \Omega _{ell}$.
Since $\Omega _{ell}$ is $G_0$-invariant
now $G_0.\exp(i{\mathfrak a_0}).x_0\subset \Omega _{ell}$.
The opposite inclusion follows via classical methods.
We have proved

\begin {proposition} \label {elliptic}
$\Omega _{ell}=G_0.\exp(i\ga_0.x_0)$
\end {proposition}

The following is a key ingredient for understanding the $G_0$--orbit
structure in bd$(\Omega _{AG})$.

\begin {proposition} \label {convex}
$\exp(i{\mathfrak a_0}).x_0\cap \text {c}\ell(\Omega _{AG})=
\text{c}\ell (\exp(i\omega _0).x_0)$
\end {proposition}

\noindent {\bf Proof.}
If $x\in \text{c}\ell (exp(i\omega _0).x_0)$, then it is elliptic
and therefore its orbit $G_0.x$ is closed. In other words
$\exp(i{\mathfrak a_0}).x_0\cap \text{c}\ell (\Omega _{AG})\supset
\text{c}\ell (\exp(i\omega _0).x_0)$.

For the opposite inclusion, observe that if
$s,s'\in \exp(i{\mathfrak a_0}).x_0$ and $s'\in G_0.s$, then
$s'=k_0(s)$ for some element $k_0$ of the Weyl group.  Thus,
if $s\in \text{c}\ell (\exp(i\omega _0).x_0)$, then
$s'\in \text{c}\ell (\exp(i\omega _0).x_0)$ as well.  Therefore, in
order to prove the opposite inclusion it is enough to show that,
given $s'\in \exp(i{\mathfrak a_0}.x_0)\cap \text{c}\ell (\Omega _{AG})$,
there exists $s\in \text{c}\ell(\exp(i\omega _0).x_0)$ with
$s'\in G_0.s$.

Given $s'$ as above, there exist sequences
$\{ s_n\}\subset \exp(i{\omega _0}).x_0$ and
$\{ s'_n\}\subset \Omega _{AG}$ such that $s'_n\in G_0.s_n$,
$s'_n\to s'$ and $s_n\to s\in \text{c}\ell (\exp(i{\omega _0}).x_0)$.
Consider the (real) categorical quotient map
$\pi : \text{Aut}_\R({\gg})\to  \text{Aut}_\R({\gg})/$ \hskip -6 pt $/ G_0$\,.
It is continuous, the base is Hausdorff and in every fiber there
is exactly one closed $G_0$-orbit.  Since $\pi (s_n)=\pi (s'_n)$,
it follows that $G_0.s=G_0.s'$.
\hfill $\square$

\begin {corollary}
Let $\Omega _{\text{c}\ell}$ denote $\{x\in \Omega :G_0.x \text { is closed }\}$.
Then \hfill \newline
\centerline{$\Omega _{\text{c}\ell }\cap \text{c}\ell (\Omega _{AG})$ =
$G_0.\text{c}\ell (\exp(i\omega _0).x_0)$ =
$\Omega_{ell}\cap \text{c}\ell (\Omega _{AG})$.}
\end {corollary}

\noindent {\bf Proof.}
From Proposition \ref{elliptic}, $\Omega _{AG}\subset \Omega _{ell}$.
By continuity, the semi-simple part of $\eta (x)$ is elliptic for
every $x\in \text{c}\ell (\Omega _{AG})$.  Thus
$\Omega _{\text{c}\ell }\cap \text{c}\ell (\Omega _{AG})\subset
\Omega_{ell}\cap \text{c}\ell (\Omega _{AG})$,
because elements of closed orbits are semi-simple.
Proposition \ref {elliptic} gives
$\Omega _{\text{c}\ell }\cap \text{c}\ell (\Omega _{AG})\subset
G_0.\text{c}\ell (\exp(i\omega _0).x_0)$, and
from Proposition \ref {convex} it follows that 
$G_0.\text{c}\ell (exp(i\omega _0).x_0)\subset
\Omega _{\text{c}\ell }\cap \text{c}\ell (\Omega _{AG})$.
So we have
$\Omega _{\text{c}\ell }\cap \text{c}\ell (\Omega _{AG})\subset
G_0.\text{c}\ell (\exp(i\omega _0).x_0)=
\Omega_{ell}\cap \text{c}\ell (\Omega _{AG}).$
Finally, if $x\in \Omega _{ell}$, then in particular it is
semi-simple and $G.x$ is closed. This proves the remaining
inclusion.  \hfill $\square$
\medskip

\noindent {\bf Existence of a $Q_2$-slice.}
\medskip

{\sc Hilbert Lemma.} 
We have fixed ${\mathfrak a_0}$\,.  Let $G_0.y \subset \text{bd}(\Omega _{AG})$
be a non--closed orbit.  We determine
a base point $z\in G_0.y$ with $\eta (z)=u\cdot s$ such that the
point $x_1$ corresponding to $s$ is in bd$(\exp(i\omega _0).x_0)$.
That uses the nilpotent element $\nu$ associated to $u$, and then
an ${\mathfrak sl}_2$-triple $(\nu = e,h,f)$ defined over $\mathbb R$
with $z$ in the orbit $S.x_1$ of the corresponding complex group.
We regard this as an analogue of the Hilbert Lemma in the case of actions
of reductive complex Lie groups.

\begin {lemma}  \label {Hilbert Lemma}
Every non-closed $G_0$--orbit in
$\text{\rm bd}(\Omega _{AG})$ contains a point
$z= \exp(\nu)\cdot \exp(i\xi) \cdot x_0$ where $\xi \in \text{\rm bd}(\omega _0)
\subset \ga_0$ and where $\nu \in \gg^{\eta(\exp(i\xi))}\cap i{\mathfrak g_0}$
is nonzero and nilpotent.
\end{lemma}

\noindent {\bf Proof.}
Let $\eta(y)=su$ be the Jordan decomposition and let
$\nu\in {\mathfrak g}^s\cap i{\mathfrak g}_0$ as in
Proposition \ref {lifting}.   Then
$\eta(y)= \eta(\exp(-\frac{1}{2}\nu)\exp(\frac{1}{2}\nu)\cdot y)=
\text{Ad}(\exp \nu)\circ\eta(\exp(\frac{1}{2}\nu)\cdot y)=u\cdot s.$
By Proposition \ref {closed} and Proposition \ref {elliptic} the
semisimple element $ \eta(\exp(\frac{1}{2}\nu)\cdot y)$ is elliptic.
Hence, Proposition \ref{elliptic}
implies the existence of  $g\in G$ and $\xi \in {\rm bd}(\omega _0)$
such that $\exp(\frac{1}{2}\nu)\cdot y=g^{-1}\exp(i\xi)\cdot x_0$\,.
Define 
$e:= \text{Ad}(g)(-\frac{1}{2}\nu)$.  Then
$g\cdot y=\exp (e) \exp(i\xi)\cdot y=\exp (e) \exp(i\xi)\cdot x_0.$
Finally, $e\in {\mathfrak g}^{g.s}={\mathfrak g}^{\eta(\exp(i\xi))}$,
and Lemma \ref{Hilbert Lemma} is proved.
\hfill $\square$
\medskip

Let $e$ be as above and $S$ be the complex group determined by the above
${\mathfrak sl}_2$-triple.  Direct
calculation shows that $S^0 _{x_1}\cong \mathbb C^{\,*}$.  Thus the
orbit $S.x_1$ is equivariantly biholomorphic to either the
$2$-dimensional affine quadric or to the complement of a smooth
quadric curve in $\mathbb C\mathbb P_2$.
The former is naturally realized as the complement of the diagonal
in $\mathbb C\mathbb P_1\times \mathbb C\mathbb P_1$ and the latter
as its quotient by the $\mathbb Z _2$-action defined by reversing
its factors.  At the level of homogeneous spaces the former is
$S/T$, where $T$ is the complexification of a (compact) maximal
torus $T_0$, and the latter is $S/N$, where $N$ is the normalizer
of $T$ in $S$.
\smallskip

In order to remind the reader of the connection to the quadric,
the orbit $S.x_1$ is referred to as a $Q_2$-slice whenever
its intersection with $\Omega _{AG}$ is an $\Omega _{AG}$ associated
to the real form $S_0$ of $S$ which is defined by the restriction of
$\sigma $, the complex conjugation of $G$ over $G_0$\,.
\medskip

{\sc Genericity.} 
Starting with an arbitrary $G_0$--orbit on bd$(\Omega _{AG})$
it is difficult to determine how the $S$-orbit $S.x_1$ intersects
$\Omega _{AG}$\,.  However, for $G_0$--orbits of {\it generic} boundary 
points this will be relatively straightforward.
Here a point $z\in \text{bd}(\Omega _{AG})$ is called {\it generic}
if $G_0.z$ is not closed and if the point $x_1 = \exp(i\xi)x_0$
constructed above by the Hilbert Lemma is a
smooth point of $\text{bd}(\exp(i\omega _0).x_0)$.
\smallskip

Recall that $\text{bd}(\omega _0)$ is defined by hyperplanes.
We refer to the points in $\text{bd}(exp(i\omega _0).x_0)$ which
correspond to points which are contained in two or more hyperplanes
as {\em corners}.  Let ${\cal E}$ denote the $G_0$-saturation
of the set of such corners, i.e., points $z$ of $\text{bd}(\Omega _{AG})$
such that $G_0.z$ has a corner in its closure.  Let
${\cal C}$ be the set of points $z\in \text{bd}(\Omega _{AG})$ such
that $G_0.z$ is closed, so $\cC = G_0.\text{bd}(exp(i\omega _0).x_0)$.
Finally, let $\text{bd}(\Omega _{AG})_{gen}$ be the closure of the complement
of ${\cal E}\cup {\cal C}$ in $\text{bd}(\Omega _{AG})$.  It follows
that this is contained in the set of generic points in the above sense.
A careful look at the invariant theory for the $G_0$-action,
leads to the following density result.

\begin {proposition}
The set $\text{\rm bd}(\Omega _{AG})_{gen}$ is open and dense in
$\text{\rm bd}(\Omega _{AG})$.
\end {proposition}

Using detailed knowledge of the $G_0$-isotropy groups at the
smooth points of $\text{bd}(\exp(i\omega _0).x_0)$, one proves
the desired result on the existence of $Q_2$-slices.

\begin {proposition}
At every generic boundary point there exists a $Q_2$-slice.
\end {proposition}

Given $y\in \text{bd}(\Omega _{AG})_{gen}$ we of course move
it via $G_0$ to an optimal point $z$ in $G_0.y$ so that $S.z$ contains
the smooth boundary point $x_1\in \text{bd}(\exp(i\omega _0).x_0)$
as above.
\vfill\pagebreak

\noindent {\bf Analysis of a $Q_2$-slice.} 
\medskip

Let $S=SL_2({\mathbb C})$ and let $S_0=SL_2({\mathbb R})$ be
embedded in $S$ as the subgroup of real matrices. As usual 
$K_0 =SO_2(\mathbb R)$ and $K = SO_2(\C)$.  
$D_0$ and $D_\infty$ denote the open $S_0$-orbits in
$\mathbb C\mathbb P_1$ containing the respective $K_0$--fixed points 
$0$ and $\infty$.
\smallskip

$S$ acts diagonally on
$Z=\mathbb C\mathbb P^1\times \mathbb C\mathbb P^1$ with one
open orbit $\Omega $, the complement of the
diagonal diag$(\mathbb C\mathbb P^1)$ in $Z$.  $\Omega $ is the complex
symmetric space $S/K$.
Note that in $\mathbb C\mathbb P^1\times \mathbb C\mathbb P^1$
there are 4 open $SL_2(\mathbb R)\times SL_2(\mathbb R)$--orbits:
the bi-disks $D_\alpha \times D_\beta$ for any pair
$(\alpha,\beta)$ from $\{0,\infty\}$.
As $S_0$-spaces, the domains $D_0\times D_\infty$ and $D_\infty\times D_0$
are equivariantly biholomorphic; further, they are subsets
of $\Omega $, and the Riemannian symmetric space $S_0/K_0$ sits in
each of them as the totally real orbit
$S_0(0,\infty)$ or $S_0(\infty,0)$.
Each can be considered as the
$\Omega_{AG}$ associated to $S$ and the real form $S_0$:
$$
\Omega _{AG}=D_0\times D_\infty=S_0\cdot \exp i\omega _0\cdot
(0,\infty)\qquad \text{or } \qquad
D_\infty\times D_0=S_0\cdot \exp i\omega _0\cdot (\infty,0).
$$
We choose the first:
 $\omega _0=(-\frac{\pi}{4},\frac{\pi}{4})h_\alpha$ where
$h_\alpha\in {\mathfrak a}$ is the normalized
coroot, i.e., $\alpha(h_\alpha)=2$.
\medskip

Our main point here is to understand $S_0$-invariant Stein domains
in $\Omega $ which properly contain $\Omega _{AG}$.
By symmetry we may assume
that said domain has non-empty intersection with $D_0\times D_0$.
Observe that
$(D_0\times D_0) \cap \Omega=D_0\times D_0\smallsetminus \text{diag}(D_0)$,
and other than $\text{diag}(D_0)$ all $S_0$--orbits in
$D_0\times D_0$ are closed real hypersurfaces.   If
$p \in D_0\times D_0\smallsetminus \text{diag}(D_0)$ let $\Omega (p)$ be the domain
bounded by $S_0\cdot p$ and diag$(D_0)$. We shall show that a function
which is holomorphic in a neighborhood of $S_0\cdot p$
extends holomorphically to $\Omega (p)$.
\smallskip

$\Sigma :=\{ (-s,s):0\le s<1\} \subset D_0\times D_0$
is a geometric slice for the action of $S_0$ on $D_0\times D_0$\,.
We say that a ($1$-dimensional) complex curve
$C\subset {\mathbb C}^2 \subset Z$ is a {\em supporting curve} for
${\rm bd }(\Omega(p))$ at $p$ if
$C\cap \text{\rm c}\ell(\Omega(p))=\{ p\} $.
Here $\text{\rm c}\ell(\Omega(p))$ denotes the topological closure in
$D_0\times D_0$\,.

\begin{proposition}
If $p \in D_0\times D_0\smallsetminus \text{\rm diag}(D_0)$ there is
a supporting curve for ${\rm bd}(\Omega(p))$ at $p$.
\end{proposition}

\noindent {\bf Proof.}
We consider $D_0 \subset \mathbb C$ as the unit disc.
We need only construct a supporting curve $C\subset \mathbb C^2$
at each point $p_s=(-s,s)\in \Sigma $, $s\ne 0$.  Define
$C_s:=\{ (-s+z,s+z):z\in {\mathbb C}\} $.  To prove
$C_s\cap \text{\rm c}\ell(\Omega(p_s))=\{ p_s\}$
let $d: D_0\times D_0\to \mathbb R$ be the distance function
of the Poincar\' e metric of $D_0$\,.  It is an $S_0$--invariant, and its values 
parameterize the $S_0$-orbits on $D_0\times D_0$\,.
\smallskip

We now claim that $d(-s+z,s+z)\geqq d(-s,s)=d(p_s)$ for $z\in \mathbb C$
and $(-s+z,s+z)\in D_0\times D_0$,
with equality only for $z=0$, i.e., $C_s$ touches
$\text{\rm c}\ell( \Omega(p_s))$ only at $p_s$.
To prove this, we compare
the Poincar\' e length of the Euclidean
segment seg$(z-s,z+s)$ in $D_0$ with that of seg$(-s,s)$.
Writing the corresponding integral for the length, it is clear from a
glance at the integrand that $d(-s+z,s+z)>d(-s,s)$ for $z \ne 0$.
This completes the proof.  \hfill $\square$
\medskip

From the above construction the boundary hypersurfaces $S(p)$ are strongly 
pseudoconvex from the viewpoint of diag$(D_0)$.  The
smallest Stein domain containing a $S_0$--invariant neighborhood of
$S_0(p)$ is $\Omega(p)\smallsetminus \text{diag}(D_0)$, so
the following is immediate.

\begin {corollary} \label{extend}
If $p\in D_0\times D_0\smallsetminus \text{\rm diag}(D_0) $
and $f$ is holomorphic on a
neighborhood of $S_0\cdot p$ then $f$ extends holomorphically
to $\Omega(p)\smallsetminus \text{\rm diag}(D_0)  $.  
If $p\in D_\infty\times D_\infty\smallsetminus \text{\rm diag}(D_\infty) $
the analogous statement holds.
\end {corollary}

The set of generic boundary points is the union of the two $S_0$-orbits,
$\text{\rm bd}_{gen}(D_0\times D_\infty)$ = 
$({\rm bd}(D_0)\times D_\infty) \cup (D_0\times {\rm bd}(D_\infty))$.
Let $z\in {\rm bd}(D_0)\times D_\infty$
(or $z\in D_0\times {\rm bd}(D_\infty)$, respectively).

\begin {corollary} \label{expand}
Let $\widehat \Omega \subset Q_2\subset \mathbb C\mathbb P^1
\times\mathbb C\mathbb P^1 $ be an
$S_0$-invariant Stein domain  that
contains $D_0\times D_\infty$ and its boundary point $z$.
Then $\widehat \Omega$ also  contains
$(D_0\times \mathbb C\mathbb P^1) \smallsetminus
\text{\rm diag}(\mathbb C\mathbb P^1) $
{\rm (}or also contains $(\mathbb C\mathbb P^1\times D_\infty) \smallsetminus
\text{\rm diag}(\mathbb C\mathbb P^1) $, respectively{\rm )}.
\end {corollary}

\noindent {\bf Proof.}
Let $B$ be a ball around $z$ which is contained in
$\widehat \Omega $. For $p\in B(z)\cap (D_0 \times D_0)$
sufficiently close to $z$, $S_0\cdot q
\subset \widehat \Omega$ for all
$q\in B(z)\cap(D_0 \times D_0)$.
The result follows from Corollary \ref{extend}.
\hfill $\square$
\smallskip

If $\widehat \Omega $ is as in Corollary \ref{expand}, then
the fibers of the projection of
$\widehat \Omega\subset \mathbb C\mathbb P^1\times
\mathbb C\mathbb P^1$ to the first factor $\C\mathbb P^1 $ can
be regarded as non-constant holomorphic curves
$f:\mathbb C\to \widehat \Omega $.  In particular,

\begin {corollary} \label {nothyper}
If $\widehat \Omega$ is as in {\rm Corollary \ref{expand}}, then $\widehat \Omega$
is not Brody hyperbolic.
\end {corollary}

\noindent {\bf Characterization of cycle domains.} 
\medskip

In the previous sections, in order to avoid unnecessary notation,
we have often replaced

\centerline{$\Omega _W(D)$: component of $C_0$ in 
$\{ gC_0 \mid g \in G \text{ and } gC_0 \subset D \} \subset G/\widetilde{K}$}

\noindent by its finite cover

\centerline{$\widetilde{\Omega _W(D)}$: component of $1K$ in
$\{ gK \mid g \in G \text{ and } gC_0 \subset D \} \subset G/K$.}

\noindent  Proposition \ref{cover_trivial} below, shows that in fact the
covering $G/K \to G/\widetilde{K}$ is bijective on $\widetilde{\Omega _W(D)}$,
inducing a holomorphic diffeomorphism $\widetilde{\Omega _W(D)} \to \Omega _W(D)$
and justifying the above--mentioned replacement.  For our main theorem,
however, even though we have not yet come to the proof of Proposition 
\ref{cover_trivial}, we still
consider $\Omega _W(D)$ as sitting up in $G/K$.

\begin {theorem} \label{equalities_finally}
Either we are in the hermitian holomorphic case and $\Omega _W(D)$ is ${\cal B}$ or
$\overline{\cB}$, or $\Omega _{AG}=\Omega _I=\Omega _W(D)=\Omega _S(D)$.
\end {theorem}

\noindent {\bf Proof.}
If $G_0$ is of hermitian type, then the result is contained in
Theorem \ref{herm_case}, or see \cite {HuW3} or \cite {WZ3}. 
Otherwise we have the inclusions and equalities (\ref{incl}) from 
Summary \ref{summary1}.
By Theorem \ref{main} above, $\Omega _H=\Omega _{AG}$, and consequently
all of those inclusions are equalities.
\hfill $\square$
\medskip

Finally we view the cycle space as it really is,
and verify that the standard projection $\pi :G/K\to G/\tilde K$ 
restricts to a holomorphic diffeomorphism of $\widetilde{\Omega _W(D)}$
onto $\Omega _W(D)$.  The projection $\pi$ is given as identification 
under the right action of the finite group $\Gamma = K \backslash \widetilde{K}$
on $G/K$.  $\Gamma$ permutes the components of 
$\{ gK \mid g \in G \text{ and } gC_0 \subset D \}$ = 
$\pi^{-1}(\{ gC_0 \mid g \in G \text{ and } gC_0 \subset D \})$.
Thus $\Omega _W(D)$ is the quotient of $\widetilde{\Omega _W(D)}$
by its stabilizer in $\Gamma $.
\smallskip

By Theorem \ref{equalities_finally}, $\widetilde{\Omega _W(D)} = \Omega _{AG}$.
Since $\Omega _{AG}$ is a cell, a finite group of diffeomorphisms
can act freely on it only if it is trivial.  Thus one may indeed
regard $\Omega _W(D)$ as being in $\Omega =G/K$.
\smallskip

We summarize this as follows.

\begin {proposition} \label{cover_trivial}
The restriction $\pi : \widetilde{\Omega _W(D)} \to \Omega _W(D)$ of the
projection $G/K\to G/\tilde K$ is biholomorphic.  In particular,
$\Omega _W(D)$ is a cell.
\end {proposition}

The following Corollary is contained in \cite{HuW3}.  Also see \cite{HuW4}.

\begin{corollary} \label{contractible_stein} In all cases,
$\Omega _W(D)$ is a contractible Stein manifold.
\end{corollary}
\vfill\pagebreak

\centerline{\Large \bf \underline{Part III: Applications and Open Problems.}}
\bigskip

In this Part we go more closely into applications of the complex geometric
methods described and developed in Part II.

\section{Recent Results on the Double Fibration Transform.} \label{sec11}
\setcounter{equation}{0}

We continue the discussion of double fibration transforms from Section
\ref{sec4}, taking advantage of the material just described in Part II.
\medskip

As explained above in Proposition \ref{Schubert slice}, there is an Iwasawa 
decomposition $G_0 = K_0A_0N_0$ such that the Schubert slice 
$\Sigma := A_0N_0(z_0) \subset D$ meets every cycle $C \in \Omega_W(D)$ 
transversally in a single point within $D$.  That gives a map
\begin{equation} \label{fiber_cspace}
\phi : \Omega_W(D) \to \Sigma := A_0N_0(z_0) \text{ by } 
C \mapsto (\Sigma \cap C) \in D.
\end{equation}
Note that $\phi^{-1}(z)$ consists of all cycles
$C \in \Omega_W(D)$ that contain $z$, so $\phi^{-1}(z) =
\mu^{-1}(z) =: F$.
\smallskip

Let $J_0$ denote the isotropy subgroup of $A_0N_0$ at $z_0$\,,
and let $F_0 = \mu^{-1}(z_0)$.  Note that $J_0$ acts on $F_0$\,.
Realize (\ref{fiber_cspace}) as the $A_0N_0$--homogeneous
fiber bundle $(A_0N_0) \times_{J_0} F_0 \to A_0N_0/J_0$\,.   
The subsets $\Omega_W(D)$ and 
$\widetilde{F}_0 := \{ C \in \Omega \mid z_0 \in C \}$ are
semialgebraic in $\Omega$, so their intersection $F_0$ has only
finitely many topological components.   As $\Omega_W(D)$ is simply
connected, it follows that $\Sigma =  A_0N_0/J_0$ is a solvmanifold with finite
fundamental group.  Thus $\Sigma$ is $C^\omega$--diffeomorphic to 
a cell.  Hence the fibration (\ref{fiber_cspace}) is trivial.  By Proposition
\ref{cover_trivial}, its total space is a cell.  So now the base and
total space of (\ref{fiber_cspace}) are cohomologicaly trivial, and thus
the same holds for the fiber $F$. We have proved

\begin{theorem} \label{holo_fiber_cell} {\rm \cite{HuW4}}
Let $F$ denote the fiber of the holomorphic fibration
$\mu : \cI(D) \to D$.  Then $F$ is connected and 
$H^r(F;\C) = 0$ for all $r > 0$.
In particular {\rm (\ref{buchdahl_conditions})} is
satisfied for every $q$, and the double fibration transforms 
$P :H^q(D;\cE) \to H^0(M; \cR^q(\mu^*\cE))$ are injective for all
sufficiently negative $\E \to D$.
\end{theorem}

\begin{remark} \label{extend_fiber_cell} {\rm \cite{HuW4}} {\em
The fiber space projection $\phi : \Omega_W(D) \to \Sigma$ is the 
restriction to open subsets of a holomorphic bundle projection 
$\widetilde{\phi} : AN \times_J \widetilde{F}_0 \to \cO$, as follows.
Let $\widetilde{F}_0 := \{ C \in \Omega \mid z_0 \in C \}$.
The complex submanifold $\cO = B(z_0) \subset Z$ where $B$ is a Borel subgroup
of $G$ that contains $A_0N_0$\,.  Thus $\Sigma = \cO \cap D$ is open in
$\cO$ by the discussion of Schubert cells and Schubert slices in 
Theorem \ref{triality} above.
$A$ and $N$ are the respective complexifications of $A_0$ and $N_0$\,.
$\cI(D)$ is open in $AN(\widetilde{F}_0)$, which is the total space of
$\widetilde{\phi} : AN \times_J \widetilde{F}_0 \to \cO$.  The connection
with $\phi : \Omega_W(D) \to \Sigma$ is that $\Omega_W(D) = A_0N_0(F_0)$,
which is open in $AN(\widetilde{F}_0 = AN\times_J \widetilde{F}_0)$,
where $J$ is the isotropy subgroup of $AN$ at $z_0$\,.  Since it is the restriction
of $\widetilde{\phi}$, the map $\phi : \Omega_W(D) \to \Sigma$ is
holomorphic.}
\end{remark}

Now we have adequately addressed the injectivity requirement 
(\ref{want_inj}) for the double fibration transform of a flag domain,
and we turn to the question (\ref{want_image}) of its image.
Since the Stein manifold $\Omega_W(D)$ is contractible, every
holomorphic vector bundle $\widetilde{\E} \to \Omega_W(D)$ is
holomorphically trivial, and in particular the Leray derived bundles
$\E^\dagger = \H^q(C, \mu^*(\E)|_{\nu^{-1}(C)})$ over $\Omega_W(D)$ are 
holomorphically trivial.
Here we have two requirements for (\ref{want_image}): we need
\begin{gather}
\text{a canonical choice of holomorphic trivialization of } 
\E^\dagger \to \Omega_W(D), \text{ and }
\label{canonical_holo_triv}\\
\text{an explicit (in that trivialization) system of PDE that specifies
the the image of } P. \label{image_PDE}
\end{gather}
This is work in progress.

\section{Unitary Representations of Real Reductive Lie Groups.} \label{sec12}
\setcounter{equation}{0}

In this Section we look at some of the implications of the double fibration
transform for representations of real reductive Lie groups.

Harish--Chandra's analysis of the holomorphic discrete series can be viewed
from the perspective of the double fibration transform as follows.  Let
$G_0$ be of hermitian type, $\cB = G_0/K_0$\,.  In this case, of course,
$D = \cB = \cI(D) = \Omega_W(D)$, the double fibration transform is the 
identity, (\ref{canonical_holo_triv}) is completely standard, and the
system (\ref{image_PDE}) consists of the $\overline{\partial}$ operator.  Let
$\E_\lambda \to \cB$ denote the homogeneous holomorphic hermitian vector bundle 
associated to the representation $E_\lambda$ of $K_0$ of highest weight
$\lambda$.  By use of his system of strongly orthogonal noncompact positive
roots, and the explicit holomorphic trivialization of $\E_\lambda \to \cB$,
he proves (i) a holomorphic section of $\E_\lambda \to \cB$ is $L^2(\cB)$
if and only if its $K_0$--isotypic components are $L^2(\cB)$, (ii) if some
nonzero $K_0$--isotypic holomorphic section of $\E_\lambda \to \cB$ is 
$L^2(\cB)$ then the constant section $f_\lambda$\,, value equal to the 
highest weight vector of $E_\lambda$\,, is $L^2(\cB)$, and (iii) $f_\lambda$
is $L^2(\cB)$ if and only if 
$\langle \lambda + \rho, \beta \rangle < 0$ where $\rho$ is half the sum of
the positive roots and $\beta$ is the maximal root.

Narasimhan and Okamoto \cite{NO} extended the Harish--Chandra construction
to ``almost all'' discrete series representations of a real group $G_0$ of
hermitian type, again always working over $D = \cB = \cI(D) = \Omega_W(D)$
where the double fibration transform is more or less invisible.  

The double fibration transform first became visible, at least in degenerate
form,  in Schmid's holomorphic construction of the discrete series 
(\cite{S3}, \cite{S5}).  There $Z = G/B$ for some Borel
subgroup $B$ and $D = G_0/T_0$ where $T_0$ is a compact Cartan subgroup,
$T_0 \subset K_0 \subset G_0$\,.  Only the ``real form'' $\phi : D \to G_0/K_0$ 
of the double fibration appears: there $G_0/K_0$ appears instead of the cycle
$\Omega_W(D)$; correspondingly $D$ appears instead of the incidence space
$\cI(D)$.  Injectivity of this real double fibration transform
$P_\R : H^q(D;\E) \to H^0(G_0/K_0; \E^\dagger)$ is given by Schmid's
``Identity Theorem''.  That theorem says that, under appropriate restrictions, 
a Dolbeault class $[\omega] \in H^q(D;\E)$ is zero if and only if every 
restriction $\omega$ is cohomologous to zero on every fiber of
$\phi : D \to G_0/K_0$\,.  This was extended a bit by Wolf \cite{W3}, for 
flag domains of the form $D \cong G_0/L_0$ with $G_0$ general reductive 
and $L_0$ compactly embedded in $G_0$\,.

The double fibration transform first appeared in modern form in the paper
\cite{WeW} of Wells and Wolf on Poincar\' e series and automorphic cohomology.
The only restriction there was that $D \cong G_0/L_0$ with $L_0$ compact,
and a small extension of the Identity Theorem was used to, in effect, prove
injectivity of the double fibration transform.  

The Penrose transform applies to the case $D = SU(2,2)/S(U(1) \times U(1,2))$.
There $L_0$ is noncompact, and perhaps that is the first such case to be
studied carefully.  See \cite{BE}.  Background work on interesting flag 
domains with
noncompact isotropy includes, of course, parts of Berger's classification
\cite{Be} of semisimple symmetric spaces, Wolf's study \cite{W1} of isotropic
pseudo--riemannian manifolds, and of course \cite{W2}.  Important cases of
construction of unitary representations using
double fibration transforms on flag domains with noncompact isotropy were 
studied in Dunne--Zierau \cite{DZ} and Patton--Rossi \cite{PR2}. This
area was first studied systematically in Wolf--Zierau \cite{WZ2}.

Finally, as noted in \cite{W6}, there are indications of a strong relation 
between the double fibration transforms of \cite{WZ2} and the construction of
unitary representations by indefinite harmonic theory of (Rawnsley, Schmid \& 
Wolf \cite{RSW}).
 
\section{Variation of Hodge Structure.} \label{sec13}
\setcounter{equation}{0}

In this Section we indicate the connection between Griffiths' theory of 
moduli spaces for compact K\" ahler manifolds (period matrix domains and 
linear deformation spaces), on the one hand, 
and flag domains, cycle spaces and double fibration transforms on the 
other hand.  Along the way we will sketch some relevant aspects automorphic
cohomology theory as developed by Wallach, Wells, Williams and Wolf.
\medskip

For Griffiths' theory see \cite{Gri1} and \cite{Gri2}.  There are expositions
contained in \cite{Gri3}, \cite{S3}, \cite{We1}, and \cite{We2}.
The first Hodge--Riemann bilinear relation specifies
a complex flag manifold $Z = G/Q$ and the second Hodge--Riemann bilinear
relation specifies an open $G_0$--orbit $D \subset Z$.
\smallskip

Let $X$ denote a compact K\" ahler manifold, $H_0^r(X;\C)$ and $H_0^r(X;\R)$
the complex and real spaces of primitive cohomology classes in degree $r$,
and $H_0^r(X;\C) = \sum_{p+q=r} H_0^{p,q}(X;\C)$ the decomposition by
bidegree.  This specifies the Hodge filtration 
$(F^0 \subset F^1 \subset \dots \subset F^r)$ of $H_0^r(X;\C)$, 
where $F^s = \sum_{i < s} H_0^{r-i,i}(X;\C)$, and thus
the complex flag $\cF(X) = (F^0 \subset F^1 \subset \dots \subset F^u)$
where $u$ is the integer part of $(r-1)/2$.
We have a nondegenerate bilinear form $b$ on $H_0^r(X;\C)$ given (on Dolbeault
representative differential forms) by 
$b(\xi,\eta) = (-1)^{r(r+1)/2} \int \omega^{n-r} \wedge \xi \wedge \eta$.
Here $\omega$ is the K\" ahler form of $X$.  Evidently
$b(H_0^{p,q}(X;\C), H_0^{p',q'}(X;\C)) = 0$ unless $p+p' = r = q+q'$.
Define $w(\xi) = (\sqrt{-1})^{p-q}\xi$ for $\xi \in H_0^{q,p}(X;\C)$.
One can formulate the Hodge--Riemann bilinear relations as (1) $b$ pairs 
$H_0^{p,q}(X;\C)$ with its complex conjugate $H_0^{q,p}(X;\C)$ and 
(2) $h(\xi,\eta) := b(w\xi, \overline{\eta})$ is positive definite on
$H_0^r(X;\C)$.
\smallskip

If $r$ is even, say $r = 2t$, then $b$ is symmetric.  It is positive definite 
on $H_0^{r-i,i}(X;\C) \oplus H_0^{i,r-i}(X;\C)$ for $i < t$, negative definite
on $H_0^{t,t}(X;\C)$.  The (identity
component of the) isometry group of $(H_0^r(X;\C),b)$ is the complex
special orthogonal group $G = SO(2h+k;\C)$ where $k = \dim  H_0^{t,t}(X;\C)$ 
and $h = \sum_{i < t} h_i$ with $h_i = \dim H_0^{r-i,i}(X;\C)$.  The
dimension sequence of the flag $\cF(X)$ specifies the complex flag manifold
$Z = G/Q$ consisting of all the flags 
$\cE = (E^0 \subset E^1 \subset \dots \subset E^{t-1})$ in $H_0^r(X;\C)$ 
with $b(E^{t-1},E^{t-1}) = 0$.  The (identity component of the) isometry group 
of $(H_0^r(X;\R),b)$ is the identity component $G_0 = SO(2h,k)^0$
of the real special orthogonal group $SO(2h,k)$.
The second bilinear relation above shows that the isotropy subgroup $L_0$
of $G_0$ at $\cF(X)$ is compact.  It follows that $L_0$ is of the
form $(U(h_0) \times \dots \times U(h_{t-1}) \times SO(k)$.
The flag $\cF(X)$ ranges (as $X$ varies) in the open $G_0$--orbit 
\begin{equation*}
D = \{\cE \mid b \ggg 0 \text{ on } E^{t-1} + \overline{E^{t-1}}\}
\cong SO(2h,k)/(U(h_0) \times \dots \times U(h_{t-1}) \times SO(k)). 
\end{equation*}
Here $U(h_i)$ preserves $(H_0^{r-i,i}(X;\C) + H_0^{i,r-i}(X;\C)) 
\cap H_0^r(X;\R)$, and $SO(k)$ preserves $ H_0^{t,t}(X;\R)$.  
\smallskip

If $r$ is odd, say $r = 2t-1$, then $b$ is antisymmetric, so $H_0^r(X;\C)$ 
has even dimension $2m$ and the isometry group of $(H_0^r(X;\C),b)$ is the
complex symplectic group $G = Sp(m;\C)$.  The dimension sequence of the 
flag $\cF(X)$ specifies the complex flag manifold $Z = G/Q$ consisting of 
all the flags
$\cE = (E^0 \subset E^1 \subset \dots \subset E^t)$ in $H_0^r(X;\C)$
with $b(E^t,E^t) = 0$.  The isometry group of $(H_0^r(X;\R),b)$ is the real
symplectic group $G_0 = Sp(m;\R)$.  As above, $G_0$ has compact isotropy
subgroup $L_0$ at $\cF(X)$, necessarily of the form 
$U(h_0) \times \dots \times U(h_t)$.  The flag $\cF(X)$ ranges (as $X$ varies) 
in the open $G_0$--orbit $D = \{\cE \mid b \text{ nondegenerate on each }
(H_0^{r-i,i}(X;\C) + H_0^{i,r-i}(X;\C))\}$, which is realized as
$Sp(m;\R)/(U(h_0) \times \dots \times U(h_t))$ where 
$h_i = \dim H_0^{r-i,i}(X;\C)$ as before.  
\smallskip

Since $G_0$ has compact isotropy subgroup $L_0$ on $D$,
we have $L_0 \subset K_0$\,, and the holomorphic double fibration 
(\ref{flag_doublefibration}) is supplemented by maps $D = G_0/L_0 \to G_0/K_0 
\subset \Omega_W(D)$.
\smallskip

Choose a basis $\{\gamma_1 , \dots , \gamma_v\}$ of the space 
$H_r(X;\Z)/\text{(torsion)}$ of $r$--cycles on $X$.  Given $\cF(X)$ we
have a basis $\{\omega^1 , \dots , \omega^u\}$ of $H_0^{r,0}(X;\C)$, 
then $H_0^{r-1,1}(X;\C)$, continuing through the $b$--isotropic space of 
the flag $\cF(X)$.  That defines a $u \times v$ period matrix
$$
\Pi(X) := \begin{pmatrix} 
  \int_{\gamma_1} \omega^1 & \dots & \int_{\gamma_v} \omega^1 \\
  \vdots                   &       & \vdots                    \\
  \int_{\gamma_1} \omega^u & \dots & \int_{\gamma_v} \omega^u
\end{pmatrix}
$$
which of course specifies $\cF(X)$.  As in
the case of period matrices of Riemann surfaces, one can change the
basis $\{\gamma_i\}$ by any integral element of $G_0$ and change the
basis $\{\omega^j\}$ by any element of $G_0$ that does not change $\cF(X)$.
Thus the moduli space for $r$--forms of compact K\" ahler manifolds $X$ with
given Hodge numbers $h_0^{p,q} := \dim H_0^{p,q}(X;\C)$, $p+q=r$, is the 
arithmetic quotient
\begin{equation} \label{moduli_space}
\begin{aligned} 
\Gamma \backslash D 
& = G_\Z \backslash G_0/L_0 \\
& = SO(2h,k;\Z)\backslash SO(2h,k)/(U(h_0) \times \dots \times U(h_{t-1}) 
	\times SO(k)) \text{ for } r \text{ even,} \\
& = Sp(m;\Z)\backslash Sp(m;\R)/(U(h_0) \times \dots \times U(h_t)) 
	\text{ for } r \text{ odd,}
\end{aligned}
\end{equation}
where $h_i = h_0^{r-i,i}$ and $\Gamma = G_\Z$ is defined by the lattice 
$H_r(X;\Z)$ in $H_r(X;\R)$.  A variation of Hodge structure of $X$ corresponds 
to a deformation of  K\" ahler structure of $X$, that is a fiber space
$\psi : U \to V$ and a distinguished point $v_0 \in V$ such that the
$X_v = \psi^{-1}(v)$ are compact K\" ahler (or algebraic) manifolds,
$X = X_{v_0}$\,, with $h_0^{p,q}(X_v) = h_0^{p,q}(X)$, and such that the
$X_v$ vary holomorphically (or algebraically).  That defines a holomorphic map
of $V \to \Gamma \backslash D$.
\medskip

Classically one constructs automorphic functions on $\Gamma \backslash D$ 
as quotients of $\Gamma$--invariant sections of holomorphic line bundles 
over $D$ (automorphic forms of a given weight).  Also classically $D$ is a 
bounded symmetric domain $Sp(g;\R)/U(g)$ and one works in a fixed holomorphic 
trivialization of the line bundles over $D$, so the $\Gamma$--invariance 
condition is expressed by a transformation law.  In this way one constructs
the function field of the moduli space $\Gamma \backslash D$.

The classical theory of automorphic functions must be modified in our context 
because in general $D$ has no nonconstant holomorphic functions \cite{W2},
and in general nontrivial homogeneous vector bundles over $D$ have no 
nonzero holomorphic sections.
Instead one considers sufficiently negative homogeneous holomorphic
vector bundles $\E \to D$.  Roughly speaking, those are the bundles whose
$L^2$ cohomology, and whose sheaf cohomology, viewed as $G_0$--modules, 
have the same underlying Harish--Chandra module.  
Their cohomology occurs in degree $\dim C_0$ where $C_0 \cong K_0/L_0$\,.
One looks for {\em automorphic cohomology}, meaning $\Gamma$--invariant 
classes in $H^q(D;\cO(\E))$.  That is a bit remote from the idea of a
function field for $\Gamma \backslash D$, but the double fibration
transform $P : H^q(D;\cO(\E)) \to H^0(\Omega_W(D);\cO(\E^\dagger))$ and the
holomorphic trivialization of $\E^\dagger \to \Omega_W(D)$ carry the
automorphic cohomology space $H^q(D;\cO(\E))^\Gamma$ to a space of
holomorphic functions $\Omega_W(D) \to H^q(C_0;\cO(\E|_{C_0}))$ with a
certain transformation law under $\Gamma$.  In this sense 
$\Gamma \backslash \Omega_W(D)$ can be a good replacement for
$\Gamma \backslash D$ as universal deformation space.

In much of the literature one considers only the situation where $G_0$
is of hermitian type and the bounded symmetric domain $\cB = G_0/K_0$
is used instead of $\Omega_W(D)$.  (Of course they are the same if $D$
is of hermitian holomorphic type.) When $G_0$ is not of hermitian type
then again $G_0/K_0$ is used instead of $\Omega_W(D)$, and
it is considered somewhat of an obstacle that $G_0/K_0$ is not a complex
manifold.  Our use of $\Gamma \backslash \Omega_W(D)$ addresses this point.

In connection with construction of automorphic cohomology, Wells \cite{We1}
showed by direct computation that $\Omega_W(D)$ is a Stein manifold in one
particular case ($r = 2$).  That result was extended in Wells--Wolf
\cite{WeW} to the more general situation of open $G_0$--orbits $D$ of the
form $G_0/L_0$ with $L_0$ compact, using a special case of the double
fibration transform together with somewhat general methods of
complex analysis (Andreotti--Grauert \cite{AnG}, Andreotti--Norguet
\cite{AN}, Docquier--Grauert \cite{DG}) associated to questions of
holomorphic convexity and the Levi problem.
The goal of \cite{WeW} was construction of automorphic cohomology as 
convergent Poincar\' e $\vartheta$--series $\vartheta_\Gamma(c) :=
\sum_{\gamma \in \Gamma} \gamma^*(c)$ where $c \in H^q(D;\cO(\E))$ is a
$K_0$--finite cohomology class.  The relevant estimates were derived from
semisimple representation theory, specifically from Hecht--Schmid \cite{HSc} 
and Schmid \cite{S4}, and the passage between $D$ and $\Omega_W(D)$.

This theory of  Poincar\' e 
$\vartheta$--series and automorphic cohomology later was developed quite a
bit.  According to \cite{W5}, if $\Gamma$ is any discrete subgroup of $G_0$\,,
$\E \to D$ is sufficiently negative and $1 \leqq p \leqq \infty$ then
every $\Gamma$--invariant $L^p(\Gamma \backslash D)$ class in $H^q(D;\cO(\E))$ 
can be realized as a Poincar\' e series $\vartheta_\Gamma(c)$ where
$c \in H^q(D;\cO(\E))$ is $L^p(D)$.  In particular this is close to the idea
of catching all of the function field.  The ``sufficiently'' part of the 
``sufficiently negative'' condition on $\E \to D$ is relaxed in Wallach--Wolf
\cite{WaW} by construction of an appropriate reproducing kernel.  Finite
dimensionality of automorphic cohomology was  proved by Williams 
(\cite{Wi1}, \cite{Wi2}, \cite{Wi3}), using index theory of Moscovici and
Connes,  for the case where $\Gamma \backslash D$ is compact.  Despite this 
development, automorphic cohomology has not yet been effectively applied to 
variation of Hodge structure.  We expect that the new information on the 
double fibration transform, presented above, will make a difference here.

\vskip 5 cm

\centerline{\begin{tabular}{ll}
ATH: & JAW: \\
Fakult\" at f\" ur Mathematik & Department of Mathematics \\
Ruhr--Universit\" at Bochum & University of California \\
D-44780 Bochum, Germany & Berkeley, California 94720--3840, U.S.A. \\
                               &                                      \\
{\tt ahuck@cplx.ruhr-uni-bochum.de} & {\tt jawolf@math.berkeley.edu}
\end{tabular}}

\enddocument
\end